\documentclass[final,onefignum,onetabnum]{siamonline220329}

\usepackage{subcaption}
\usepackage{graphicx}
\usepackage{hyperref}
\usepackage{amsmath}
\newtheorem{conjecture}[theorem]{Conjecture}
\newtheorem{fact}[theorem]{Fact}
\DeclareMathOperator{\Ima}{Im}
\DeclareMathOperator{\dom}{dom}
\renewcommand{\d}{\mathrm{d}}

\usepackage{lipsum}
\usepackage{amsfonts}
\usepackage{graphicx}
\usepackage{epstopdf}
\usepackage{algorithmic}
\ifpdf
  \DeclareGraphicsExtensions{.eps,.pdf,.png,.jpg}
\else
  \DeclareGraphicsExtensions{.eps}
\fi

\usepackage{enumitem}
\setlist[enumerate]{leftmargin=.5in}
\setlist[itemize]{leftmargin=.5in}

\newsiamremark{remark}{Remark}
\newsiamremark{hypothesis}{Hypothesis}
\crefname{hypothesis}{Hypothesis}{Hypotheses}
\newsiamthm{claim}{Claim}

\headers{Dynamical Fractal: Theory and Case Study}{Junze Yin}

\title{Dynamical Fractal: Theory and Case Study\footnote{This work was funded by Student Research Award (SRA) from Boston University Undergraduate Research Opportunity Program}}

\author{Junze Yin\thanks{Department of Mathematics \& Statistics, Boston University, Boston, MA 02215 USA 
  (\email{junze@bu.edu}).}}

\usepackage{amsopn}

\ifpdf
\hypersetup{
  pdftitle={Dynamical Fractal: Theory and Case Study},
  pdfauthor={Junze Yin}
}
\fi

\begin{document}

\maketitle

\begin{abstract}
Urbanization is a phenomenon of concern for planning and public health: projections are difficult because of policy changes and natural events, and indicators are multiple. There are previous studies of development that used fractals, but none for this specific problem, nor extrapolating the future trend. In the first part of this paper, we construct a theoretical framework for analyzing dynamic (changing) fractals and extrapolating their future trends based on their fractal dimension—a measure of the complexity of the fractal. We believe this approach holds enormous potential for applications in analyzing changing fractals in the real world, such as urban growth, cells, cancers, etc., all of which are invaluable to research. This theoretical framework may shed light on a factor overlooked in past research: the trend of how fractals change. In the second part of this paper, we apply this theoretical framework to study the urbanization of Boston. We compare several maps and measurements of fractal dimensions, produce code\footnote{\url{https://github.com/yinj66/fractal-dimension}} that reads the maps and divides the city into subsections, and ultimately graph the fractal dimension over time using both differential and difference equations. Finally, we postulate the logistic equation as a model to fit the evolution of the fractal dimension, as well as the total population derived from census data, which serves as a component for comparing dynamical fractals.
\end{abstract}



\section{Introduction}

Fractals do not have a rigorous mathematical definition \cite{Falconer2006}. They are geometric structures without a regular shape, but they are considered to be self-similar and have infinitely many details as we zoom in. Mandelbrot first proposed a formal approach to analyzing fractals for measuring the coastline of Britain \cite{2}. Fractal-like structures include waves, rocks, metals, rivers, plants, degree of urbanization, etc. \cite{1, 3, 4, 7, 10, 14}. The fractal dimension $D \in [1, 2]$ is a quantitative measurement of the complexity of a fractal, which is defined as:
\begin{align*}
    -D = \frac{\log {N}}{\log {\epsilon}},
\end{align*}
where $N$ denotes the number of measurement units, and $\epsilon$ represents the scaling factor.

Nowadays, existing research \cite{1,7,9,4,3} analyzes real-world fractals by computing the fractal dimension of these fractals at some chosen time points. Although this can reflect the characteristics of the fractals at certain time points, things in our world are constantly changing, so the fractal dimension at certain time points is not enough for us to visualize the long-term changes of these fractals. The trend, we believe, is a crucial component for understanding a real-world fractal and should not be overlooked. 

In this paper, we first construct a theoretical framework for analyzing the long-term and short-term trends of fractals. To do this, we give the formal definition of a dynamical fractal, which is a continuously and irregularly changing fractal. We use it to denote real-world fractals. We analyze it by modeling the fractal dimension of continuous time points. Although in real-world applications, it is only possible to analyze the fractal at finitely given time points, we use the difference equation and the differential equation to extrapolate the future trends to fill the gaps between the time points. Moreover, to address the issue of lacking fractal dimension data for extrapolating the long-term behavior of a certain dynamical fractal, we use the concept of curvature to compare the population growth, showing the extent of similarity of two dynamical fractals. If they are similar, then we believe that they have the same long-term trend.

Because all real-world fractals are changing, we believe that our theoretical framework can be applied to all of these real-world changing fractals. Fractal analysis is widely used in biology and medical image analysis \cite{zm04,fkoh92,mmf+17,kms17,nm23,nuc12,Chiang2015,ddmp10}, climate science \cite{ms12,mvm+15,o97,sl11}, and traffic and urbanization \cite{3,9,1,4}. Our theoretical framework for analyzing dynamic fractals can support these research works in various ways, including but not limited to: assisting \cite{ddmp10} in studying the trends of cancer development, helping \cite{o97, ms12} in analyzing the changes in cloud shapes to improve weather forecasting. In this paper, we apply our theoretical concepts to analyze the long-term and short-term trends of urbanization in Boston.

Urbanization is on the rise all over the globe, according to the National Geographic Society \cite{20}. Dahly and Adair proposed that the rapid urbanization of the growing world has a huge effect on human health \cite{19}. It is necessary to find the degree of urbanization of each city and model it mathematically to be able to predict it. Based on this prediction, city planners will have more information to make policies and plans. There are many possible indicators of urbanization, such as population growth, the number of patents, and the transportation network, some of which have been assigned a fractal dimension \cite{Bettencourt2013}. However, assigning the city a fractal dimension based on the way it looks in aerial photography is relatively new \cite{1} and required several pieces of data mining and programming techniques. It is widely held that the larger the fractal dimension, the more urbanized the city is, as shown in \cite{9, bmnc03}. Besides using the traditional methods for getting the fractal dimension data, we introduce more accurate ways to approximate the fractal dimension of a city, such as using the planimetric map, a human-made map that contains the very precise structure of the buildings, and using deep learning techniques to process the satellite images. We use these methods to approximate the fractal dimension of Boston from 2000 to 2020 and finally compare the advantages and disadvantages of each of them. Since we do not have enough fractal dimension data to measure the long-term trend of Boston, we use the data from another similar region - Manhattan. We also use the technique developed from the theoretical framework to show that the long-term trend of the fractal dimension of Manhattan is similar to that of Boston.

\paragraph{Roadmap}
In Section~\ref{sc:df}, we formally define dynamical fractal and analyze it. In Section~\ref{sec:case_study}, we start introducing our case study: the techniques for collecting data and analyzing the advantages and the disadvantages of each of these techniques. In Section~\ref{sc:atfdb}, we apply our theoretical concepts to analyze the fractal dimension of Boston, which is used as an indicator of its urbanization trend. In Section~\ref{sec:conclusion}, we make a conclusion and discuss future research areas. In Appendix~\ref{sec:preli}, we introduce the background of this paper: the basic mathematical definitions and properties used to support our theoretical framework.

\section{Dynamical Fractal}
\label{sc:df}

In Section~\ref{sc:df:fractal}, we introduce the basic properties of fractals. In Section~\ref{sc:df:dynamical_fractal}, we give the formal definition of the dynamical fractal and analyze its properties. In Section~\ref{sc:df:model}, we present the techniques of interpreting the trend of fractal dimension, an indicator measuring the change of the dynamical fractal, both in the short term and in the long term. In Section~\ref{sc:dccef}, we use curvature as a tool to compare the population growth of different dynamical fractals to learn about their long-term behaviors. 

\paragraph{Notations}

We use $\mathbb{N} = \{0, 1, 2, \dots\}$ to denote the set of natural numbers. $\mathbb{R}$, $\mathbb{Z}$, and $\mathbb{Z}^+$ represents the set of real numbers, the set of integers, and the set of positive integers, respectively. Let $n \in \mathbb{Z}^+$. $\mathbb{R}^n$ represents the set containing all $n$-dimensional vectors whose entries are all real numbers. Let $X$ be a set. $\mathcal{P}(X)$ represents the power set of $X$, namely $\mathcal{P}(X) = \{x \mid x \subseteq X\}$. $X^c$ represents the complement of $X$. A function $f$ is of $C^\infty$ if $f$ is infinitely continuous and differentiable. We use $\frac{\d f}{\d t}$ or $\dot{f}$ to represent the derivative of $f$ with respect to $t$ and use $\ddot{f}$ to represent the second derivative of $f$ with respect to $t$. We use $\dom(f)$ and $\Ima(f)$ to denote the domain and the image of $f$, respectively. For a vector $v \in \mathbb{R}^n$, we use $v_i$ to denote its $i$-th entry and use $\|v\|$ to denote the $\ell_2$ norm of $v$, namely $\|v\| := \sqrt{\sum_{i = 1}^n v_i^2}$. For $v_1, v_2 \in \mathbb{R}^3$, we use $\times$ to denote the cross product, where $\|v_1 \times v_2\| := \|v_1\| \cdot \|v_2\| \cdot \sin(\theta)$ and $\theta$ is the angle between $v_1$ and $v_2$.

\subsection{Fractals}
\label{sc:df:fractal}

In this section, we introduce the basic concepts of fractals, which serve as the foundation of our theoretical analysis in this paper.

Fractals are irregular geometric structures that cannot be described by using classical geometry, and each part of which has the same statistical character as the whole, which is defined as self-similarity \cite{2}. To better analyze them, the fractal dimension is introduced as an indicator of the complexity of the fractals. It is able to reflect whether a certain structure is ``dense or scattered, compact or dispersed" \cite{1}. Also, it has the same property as the dimension that  
\begin{align}
    N = \epsilon^{-D},
    \label{eq1}
\end{align}
where $D\in[1,2]$, $N$ is the number of pieces that the geometric object is evenly separated into, and $\epsilon$ is the scaling factor. For example, if a square, whose side length is 1, is separated into 4 pieces, each of which has the side length $1/2$, which can be expressed by $4 = (\frac{1}{2})^{-2}$. By placing $D$ on the left side of this equation, we can get:
\begin{align}
    -D = \frac{\log {N}}{\log {\epsilon}}.
    \label{eq:2.2}
\end{align}

However, when analyzing real-world fractals, it is not possible to get the exact value of $N$ and $\epsilon$. Instead, two fractal dimension approximation methods have been used: box-counting analysis and radial analysis\footnote{There are two more methods called dilation analysis and correlation analysis which are not discussed in this paper.} \cite{14}. We will introduce and use these methods in Section \ref{sec:case_study} for our case study.  

\subsection{Dynamical Fractal and Its Properties}
\label{sc:df:dynamical_fractal}

In this section, we introduce the formal definition of the dynamical fractal and analyze its properties. 

Informally, we define a dynamical fractal as an irregularly and continuously changing fractal. Before providing its formal definition, we need to define what an irregularly changing fractal and a continuously changing fractal are.

\begin{definition}[(Continuously) changing fractal]\label{def:changing_fractal}
    Let $f$ be a fractal. 

    Let $\phi_f : (t_0, t_*) \to \mathbb{R}$ be a function, where $\phi_f(t)$ represents the fractal dimension of the fractal $f$ at time $t \in (t_0, t_*) \subseteq \mathbb{R}$.

    $f$ is said to be a changing fractal on $(t_0, t_*)$ if there exists $t_1, t_2 \in (t_0, t_*)$ with $t_1 \neq t_2$ such that
    \begin{align*}
        \phi_f(t_1) \neq \phi_f(t_2).
    \end{align*}

    $f$ is said to be a continuously changing fractal on $(t_0, t_*)$ if $f$ is a changing fractal and $\phi_f$ is continuous on $(t_0, t_*)$.
\end{definition}

Then, we give the definition of the irregularly changing fractal.

\begin{definition}[Irregularly changing fractal]\label{def:irregularly_changing_fractal}
    Let $f$ be a changing fractal in $(t_0, t_*)$ (see Definition~\ref{def:changing_fractal}). 

    Let $\psi_f : (t_0, t_*) \to \mathbb{R}$, satisfying 
    \begin{align*}
        \psi_f(t) := \begin{cases}
            t & \text{if $f$ is not a fractal at } t\\
            t_* & \text{if $f$ is a fractal at } t
        \end{cases}
    \end{align*}

    $f$ is said to be an irregularly changing fractal in $(t_0, t_*)$ if the set
    \begin{align*}
        \Ima(\psi_f) \cap \dom(\psi_f)
    \end{align*}
    is countable. 
\end{definition}

The intuition for defining this irregularly changing fractal is to make the dynamical fractals to be fractal for almost all time. We want to avoid dynamical fractals becoming non-fractals in a certain time interval as time changes, but we allow situations in which dynamical fractals become non-fractals countably many times. 

\begin{definition}[Dynamical fractal]\label{def:dynamical_fractal}
    A dynamical fractal $d_f$ defined in $(t_0, t_*) \subseteq \mathbb{R}$ is an irregularly and continuously changing fractal in $(t_0, t_*)$ (see Definition~\ref{def:changing_fractal} and \ref{def:irregularly_changing_fractal}). Mathematically, we define it to be the set containing all the changes at any given time, namely
    \begin{align*}
        d_f := \{d_{f, t} \mid t \in (t_0, t_*)\},
    \end{align*}
    where $d_{f, t}$ is the dynamical fractal $d_f$ at the time $t$. We use $\mathcal{D}$ to denote the set containing all the dynamical fractals.
    \label{df:3.1}
\end{definition}

Now, we analyze the properties of the dynamical fractal.

\begin{conjecture}
    All fractals existing in our world are dynamical fractals.
    \label{lm:3.2}
\end{conjecture}

Everything in our world is changing continuously. Therefore, fractals in our world are continuously changing fractal (see Definition~\ref{def:changing_fractal}). We assume that their changes are irregular, then by Definition \ref{df:3.1}, they are dynamical fractals. 

This assumption is empirical. Thinking about the rocks, coastlines, and city borderlines, all of these are fractals for all the time we observe them.

\begin{lemma}\label{lem:2.9}
    For an arbitrary fractal $d_f$ existing in our living world in a given time $t \in (t, t_*)$, where $t, t_* \in \mathbb{R}$, there exists a well-defined function 
    \begin{align*}
        \Phi : d_f \to \mathbb{R},
    \end{align*}
    which maps all $d_{f, t}$ to a unique fractal dimension.
    \label{thm:3.3}
\end{lemma}

\begin{proof}
    By Conjecture~\ref{lm:3.2}, $d_f$ is a dynamical fractal. 

    Let $\Phi: d_f \to \mathbb{R}$ be defined as
    \begin{align}\label{eq:phi}
        \Phi(d_{f, t}) = \phi_f(t),
    \end{align}
    where $\phi_f$ is defined as in Definition~\ref{def:changing_fractal}.

    Suppose there exists $t_1, t_2 \in (t, t_*)$ with $t_1 = t_2$.

    Then, by Definition~\ref{def:changing_fractal} we have
    \begin{align*}
        \phi_f(t_1) = \phi_f(t_2),
    \end{align*}
    which implies that
    \begin{align*}
        \Phi(d_{f, {t_1}}) = \Phi(d_{f, {t_2}}),
    \end{align*}
    by Eq.~\eqref{eq:phi}.
\end{proof}

\begin{definition}[Self similarity of a dynamical fractal]\label{def:self_similar}
    A dynamical fractal is said to be self-similar if for almost all time $t$, the fractal is self-similar.
\end{definition}

\begin{theorem}\label{thm:self_similar}
    All dynamical fractals in our world are self-similar.
    \label{thm:3.6}
\end{theorem}

\begin{proof}
    For any given dynamical fractal $d_f$ and for almost all time $t$, $d_{f, t}$ is a fractal, by Conjecture~\ref{lm:3.2}.
    
    Since all fractals are self-similar, by Definition~\ref{def:self_similar}, $d_f$ is self-similar.
\end{proof}

\subsection{Modeling the Fractal Dimension through the Difference and Differential Equations}
\label{sc:df:model}

In this section, we introduce the methods of analyzing the short-term and long-term trends of dynamical fractals.

Although we assume all fractals in our living world are dynamical fractals (see Conjecture~\ref{lm:3.2}), in this paper, we analyze urbanization, so we only treat cities as dynamical fractals. Since the dynamical fractal is closely related to the change of fractal dimension, to study it, we focus on modeling the change of the fractal dimension. We use the difference equation to extrapolate the discrete points of fractal dimension. Therefore, we first show that in a given dynamical fractal $d_f$ defined in $(t_0, t_*)$, for all $t$ we pick from $(t_0, t_*)$, $d_{f, t}$ is a fractal.

\begin{lemma}
    Let a dynamical fractal $d_f$ be defined in $(t_0, t_*) \subseteq \mathbb{R}$.

    Let $t_1, t_2, \dots, t_n \in (t_0, t_*)$ be the time we choose.

    Then, for all $t_i \in (t_0, t_*)$, where $i \in \{1, \dots, n\}$, $d_{f, t_i}$ is a fractal.
\end{lemma}

\begin{proof}
    Let the probability space $((t_0, t_*), \Sigma, P)$ be defined in Fact~\ref{fac:probability}.

    Let $C$ be the set, where 
    \begin{align*}
        C := \{c \mid d_{f, c} \text{ is not a fractal}\}.
    \end{align*}

    By the definition of the dynamical fractal (see Definition~\ref{def:dynamical_fractal}), we have that $C$ is countable, which implies that (by Definition~\ref{def:countable}) there exists an injective function $f$ such that
    \begin{align*}
        f: C \to \mathbb{N}.
    \end{align*}

    Therefore, $c \in C$ can be indexed, namely
    \begin{align}\label{eq:C}
        C = \{t_j \mid j \in J\},
    \end{align}
    where $J$ is an index set.

    Therefore, we have
    \begin{align}\label{eq:2.4}
        P(C)
        \leq & ~ P(\bigcup_{j = 1}^\infty \{t_j\}) \notag\\
        = & ~ \sum_{j = 1}^\infty P(\{t_i\}) \notag\\
        = & ~ \sum_{j = 1}^\infty 0 \notag\\
        = & ~ 0,
    \end{align}
    where the first step follows from Eq.~\eqref{eq:C} and the index set $J \subseteq \mathbb{N}$, the second step follows from the definition of the probability space (see Definition~\ref{def:probability}), the third step follows from Fact~\ref{fac:probability}, and the last step follows from simple algebra.

    On the other hand, by the definition of probability space (see Definition~\ref{def:probability}), we have
    \begin{align}\label{eq:2.5}
        P(C) \geq 0.
    \end{align}

    Combining Eq.~\eqref{eq:2.4} and Eq.~\eqref{eq:2.5}, we have
    \begin{align*}
        P(C) = 0.
    \end{align*}

    Therefore, the probability of choosing a $c \in C$ which makes $d_{f, c}$ be a non-fractal is $0$.

    Thus, for all $t_i$ chosen, where $i \in \{1, \dots, n\}$, we have $d_{f, t_i}$ is a fractal.
\end{proof}

Now, we define the difference equation.

\begin{definition}[First order difference equation]\label{def:difference}
    Let $t_0, t_* \in \mathbb{R}$.

    Let $d_f$ be a real-world fractal, and so is a dynamical fractal (see Conjecture~\ref{lm:3.2}), which is defined in $(t_0, t_*)$.

    Let $\Phi^* : d_f \to \mathbb{R}$ be the approximation of the function $\Phi$ (Lemma~\ref{lem:2.9}). 

    Let $t_1, t_2, \dots, t_n \in (t_0, t_*)$, where $n \in \mathbb{Z}^+$.

    Then, the first-order difference equation is given by
    \begin{align*}
        \Phi^*(d_{f, t_{m + 1}}) - \Phi^*(d_{f, t_{m}}) = F(t_{m}, \Phi^*(d_{f, t_{m}})),
    \end{align*}
    where $m \in \mathbb{Z}^+$.
\end{definition}

In the long run, we expect that the fractal dimension trend of cities will be like a logistic differential equation since we believe that urbanization has a lower limit, namely the state before people build this city, and an upper limit, the state when people finish building the majority of the cities.

\begin{definition}[(Logistic) differential equations]\label{def:differential}
    Let $t_0, t_* \in \mathbb{R}$.

    Let $d_f$ be a real-world fractal, and so is a dynamical fractal (see Conjecture~\ref{lm:3.2}), which is defined in $(t_0, t_*)$.

    Let $\Phi : d_f \to \mathbb{R}$ be a function in Lemma~\ref{lem:2.9}. 

    Let $\Phi^* : d_f \to \mathbb{R}$ be the approximation of the function $\Phi$, and $\Phi^*$ is differentiable. 

    Then, the differential equation is defined as
    \begin{align*}
        \frac{\d \Phi^*}{\d t} := F(\Phi^*).
    \end{align*}

    The logistic differential equation is a particular differential equation, which is defined as 
    \begin{align*}
        \frac{\d \Phi^*}{\d t} := r \cdot \Phi^* \cdot (1 - \frac{\Phi^*}{K}),
    \end{align*}
    where $K \in \mathbb{R}$ is the carrying capacity and $r \in \mathbb{R}$.
\end{definition}

\begin{remark}
    Note that in real-world applications, we can only pick finitely many $t_i \in (t_0, t_*)$, where $i \in \mathbb{Z}^+$. 
    
    Based on the pattern of these $\Phi(d_{f, t_i})$, we find a smooth function $\Phi^*$, which is of $C^\infty$, to map the $d_{f, t_i}$ to the approximated fractal dimension. Since the purpose of this $\Phi^*$ is to extrapolate the long-term trend of the fractal dimension, making it smooth will only omit some short-term details of the fractal dimension but does not influence the long-term trend.
\end{remark} 

\begin{lemma}\label{lem:solve_logistic}
    Let 
    \begin{align*}
        \frac{\d \Phi^*}{\d t} = r \cdot \Phi^* \cdot (1 - \frac{\Phi^*}{K}),
    \end{align*}
    be a logistic differential equation defined in Definition~\ref{def:differential}.

    We define $A \in \mathbb{R}$ to be 
    \begin{align*}
        A := \frac{K - \Phi^*(d_{f, 0})}{\Phi^*(d_{f, 0})}.
    \end{align*}

    Then, $\frac{\d \Phi^*}{\d t}$ can be solved as
    \begin{align*}
        \Phi^*(d_{f, t}) = \frac{K}{1 + Ae^{-rt}}.
    \end{align*}
\end{lemma}
\begin{proof}
    The proof of this lemma is very standard, so we omit its proof. The only thing we should notice is that
    \begin{align*}
        \frac{\d \Phi^*(d_{f, t})}{\d t} 
        = & ~ \frac{\d \Phi^*(d_{f, t})}{\d d_{f, t}} \cdot \frac{\d d_{f, t}}{\d t} \\
        = & ~ \frac{\d \Phi^*(d_{f, t})}{\d d_{f, t}} \cdot 1 \\
        = & ~ \frac{\d \Phi^*(d_{f, t})}{\d d_{f, t}},
    \end{align*}
    where the first step follows from the chain rule, the second step follows from the definition of the dynamical fractal (see Definiiton~\ref{def:dynamical_fractal}), and the last step follows from simple algebra.
\end{proof}

\begin{lemma}\label{lem:difference_from_differential}
    Let the logistic differential equation 
    \begin{align*}
        \frac{\d \Phi^*}{\d t} = r \cdot \Phi^* \cdot (1 - \frac{\Phi^*}{K})
    \end{align*}
    be defined as in Definition~\ref{def:differential}.

    Let 
    \begin{align*}
        x(d_{f, t}) = \frac{r \cdot \Phi^*(d_{f, t})}{(r + 1) \cdot K}
    \end{align*}
    and $b = r + 1$.
    
    Then, we can transform the differential equation into a form of the difference equation:
    \begin{align*}
        x(d_{f, t} + 1) \approx b(1 - x(d_{f, t}))x(d_{f, t}).
    \end{align*}
\end{lemma}

\begin{proof}

Note that 
\begin{align*}
    \frac{\d \Phi^*}{\d t} 
    = & ~ \frac{\d \Phi^*}{\d d_{f, t}}\\
    = & ~ \lim_{\Delta d_{f, t} \to 0} \frac{\Phi^*(d_{f, t} + \Delta d_{f, t}) - \Phi^*(d_{f, t})}{\Delta d_{f, t}}\\
    \approx & ~ \frac{\Phi^*(d_{f, t} + 1) - \Phi^*(d_{f, t})}{1}\\
    = & ~ \Phi^*(d_{f, t} + 1) - \Phi^*(d_{f, t}),
\end{align*}
where the first step follows from the proof of Lemma~\ref{lem:solve_logistic}, the second step follows from the definition of the derivative, the third step follows from letting $\Delta d_{f, t} = 1$, and the last step follows from simple algebra.

Combining with the definition of the logistic differential equation (see Definition~\ref{def:differential}), we have
\begin{align*}
    \Phi^*(d_{f, t} + 1) - \Phi^*(d_{f, t}) \approx r \cdot \Phi^*(d_{f, t}) \cdot (1 - \frac{\Phi^*(d_{f, t})}{K}),
\end{align*}
which implies that 
\begin{align*}
    \Phi^*(d_{f, t} + 1) \approx (r - \frac{r \cdot \Phi^*(d_{f, t})}{K} + 1) \cdot \Phi^*(d_{f, t}).
\end{align*}

With 
\begin{align*}
    x(d_{f, t}) = \frac{r \cdot \Phi^*(d_{f, t})}{(r + 1) \cdot K}
\end{align*}
and $b = r + 1$, and we can get:
\begin{align*}
    x(d_{f, t} + 1) \approx b(1 - x(d_{f, t}))x(d_{f, t}).
\end{align*}

\end{proof}

To analyze the trend, we introduce the concept of stability.

\begin{definition}[Stability, in Chapter 5, Section 2D of \cite{12}]\label{def:stable}
    Let $\epsilon > 0$.

    A differential equation (see Definition~\ref{def:differential}) or a difference equation (see Definition~\ref{def:difference}) is stable on an equilibrium level $\Phi^*(d_{f, t}) = c$, where $c \in \mathbb{R}$ if there exists a neighborhood $N = (c - \epsilon, c + \epsilon)$ of $c$ such that for all $y_0 \in N$, there exists a $d_{f, t_0} \in d_f$, such that 
    \begin{align*}
        \Phi^*(d_{f, t_0})= y_0.
    \end{align*}
\end{definition}

\begin{lemma}[A revised version of Chapter 5, Section 2F of \cite{12}]\label{lem:stable}
    Let $d_{f, t}$ be a dynamical fractal (see Definition~\ref{def:dynamical_fractal}).

    Let 
    \begin{align*}
        x(d_{f, t} + 1) = b(1 - x(d_{f, t}))x(d_{f, t})
    \end{align*}
    be a difference equation (see Definition~\ref{def:difference}).

    Then, we have
    \begin{enumerate}
        \item when $b \in (1,2]$, there is a large neighborhood of stability;
        \item when $b \in (2,3)$, the stability is near the equilibrium;
        \item when $b \in (3,1+\sqrt5)$, there is no stability near the equilibrium but oscillates between two points;
        \item when $b \in [1+\sqrt5,\infty)$, the graph becomes more unstable.
    \end{enumerate}
\end{lemma}

\begin{lemma}
    Let 
    \begin{align*}
        \frac{\d \Phi^*}{\d t} = r \cdot \Phi^* \cdot (1 - \frac{\Phi^*}{K}),
    \end{align*}
    be the logistic differential equation defined in Definition~\ref{def:differential}.

    There exists a dynamical fractal $d_{f, t}$ such that
    \begin{align*}
        0 < \Phi^*(d_{f, t}) < K.
    \end{align*}

    Then, $\Phi^*$ is stable.
    \label{thm:stable}
\end{lemma}

The proof of this Lemma is very standard, so we omit it in this paper.

\begin{theorem}
    Let 
    \begin{align*}
        \frac{\d \Phi^*}{\d t} = r \cdot \Phi^* \cdot (1 - \frac{\Phi^*}{K}),
    \end{align*}
    be the logistic differential equation defined in Definition~\ref{def:differential}.

    Let 
    \begin{align*}
        x(d_{f, t} + 1) \approx b(1 - x(d_{f, t}))x(d_{f, t})
    \end{align*}
    be the corresponding difference equation of $\frac{\d \Phi^*}{\d t}$ (see Lemma~\ref{lem:difference_from_differential})

    Then, we have
    \begin{enumerate}
        \item if $r \in (0,2)$, the differential equation and its corresponding difference equation are both stable;
        \item if $r \notin (0,1]$, their stabilities are different.
    \end{enumerate}
    \label{c:stability}
\end{theorem}

\begin{proof}
    Suppose $r \in (1,2]$. By Lemma~\ref{lem:difference_from_differential}, we have $s \in (0,1]$.

    Suppose $r \in (2,3)$. By Lemma~\ref{lem:difference_from_differential}, we have $s \in (1,2)$.

    Suppose $r \in (3,1+\sqrt{5})$. By Lemma~\ref{lem:difference_from_differential}, we have $s \in (2,\sqrt{5})$.

    Suppose $r \in [1+\sqrt{5}, \infty)$. By Lemma~\ref{lem:difference_from_differential}, we have $s \in [\sqrt{5}, \infty)$.
    
    By Lemma~\ref{lem:stable}, the difference equation is stable when $b\in (1,2]$ or $b\in (2,3)$, which implies $s \in (0,2)$. 
    
    By Lemma~\ref{thm:stable}, the differential equation is always stable. 
\end{proof}

If the differential and the difference equations are both stable on the equilibrium level $\Phi^* = K$, then we expect that the fractal dimension will be convergent as $t\to\infty$. 

\subsection{Developing the Concept of Curvature to Compare the Exponential Functions}
\label{sc:dccef}

In this section, we introduce a new method that may compare the similarity between the exponential functions. Our purpose is to use that as an indicator to compare population growth. If two population growths are similar, then the fractal dimensions are expected to be similar.

In real-world applications, we might not always find enough fractal dimension data to extrapolate the long-term behavior of dynamical fractals. To address this issue, we analyze the population growth of a dynamical fractal, find another dynamical fractal that has similar population growth, and use the long-term fractal dimension trend of this dynamical fractal to represent the trend of the one that lacks the long-term fractal dimension data. 

\begin{definition}[Population]\label{def:population}
    Let $d_f$ be a dynamical fractal defined in $(t_0, t_*)$ (see Definition~\ref{def:dynamical_fractal}).

    The population of this dynamical fractal is a function
    \begin{align*}
        P : d_f \to \mathbb{R}
    \end{align*}
    where $P(d_{f, t})$ represents the population of the dynamical fractal $d_f$ at time $t$. The population is a quantitative variable and represents the amount of a certain object or organism.
\end{definition}
\begin{remark}
    Only some dynamical fractals are assigned with the population: for a certain city, the population represents the total number of people living in this city; for bacterial colonies (a dynamical fractal), the population represents the total number of bacteria. In both cases, the population may indirectly reflect the fractal dimension of the corresponding dynamical fractal: with more people living in a city, the city is considered to be more urbanized; with more bacteria living in the Petri dish, bacterial colonies are expected to be more complex. However, for dynamical fractals like rocks or coastlines, it is unlikely to find such quantitative variables which may reflect their complexity. 
\end{remark}

The population of a dynamical fractal can be modeled by exponential functions, but at different periods, the growth rates are different because of the influence of social factors. In \cite{12}, when using the exponential function to model the US population between 1900 to 1970, Meyer finds that the population can be the best fit by using two exponential functions, one covering the first three points of the population at different times and the other covering the next four points. 

\begin{definition}[Similar Population Growth]\label{def:similar}
    Let $d_f, d_g$ be dynamical fractals defined in $(t_{f0}, t_{f*})$ and $(t_{g0}, t_{g*})$, respectively (see Definition~\ref{def:dynamical_fractal}).

    Let $P_f: d_f \to \mathbb{R}$ and $P_g: d_g \to \mathbb{R}$ be the population functions for $d_f$ and $d_g$, respectively (see Definition~\ref{def:population}).

    Let $P_f^* = \{P_{f, 1}, P_{f, 2}, \dots, P_{f, n}\}$ and $P_g^* = \{P_{g, 1}, P_{g, 2}, \dots, P_{g, n}\}$ be set containing exponential functions modeling the population $P_f$ and $P_g$, respectively, where 
    \begin{align*}
        P_{f, 1} : (d_{f, t_{f0}}, d_{f, t_{f1}}) \to \mathbb{R} &\text{ and } P_{g, 1} : (d_{g, t_{g0}}, d_{g, t_{g1}}) \to \mathbb{R}\\
        P_{f, 2} : [d_{f, t_{f1}}, d_{f, t_{f2}}) \to \mathbb{R} &\text{ and } P_{g, 2} : [d_{g, t_{g1}}, d_{g, t_{g2}}) \to \mathbb{R}\\
        & ~ \vdots\\
        P_{f, n} : [d_{f, t_{f(n-1)}}, d_{f, t_{f*}}) \to \mathbb{R} &\text{ and } P_{g, n} : [d_{g, t_{g(n-1)}}, d_{g, t_{g*}}) \to \mathbb{R}.
    \end{align*}

    Their population growths of $d_f$ and $d_g$ are called similar, written as $P_f^* \sim P_g^*$, if the exponential functions in $P_f^*$ have similar shapes (see Definition~\ref{def:similar_shape}) with the exponential functions in $P_g^*$.
\end{definition}

By writing similar shapes, we mean that the ratios of the average curvatures of the exponential functions in $P_f^*$ and in $P_g^*$ are similar.

\begin{lemma}\label{lem:avergae_curvature}
    Let $P_{f, i}$ be defined in $[d_{f, t_{f(i - 1)}}, d_{f, t_{fi}})$ (see Definition~\ref{def:similar}), where $i \in \{1, 2, \dots, n - 1\}$.

    Let $\gamma_f(t) \in \mathbb{R}^3$ be defined as
    \begin{align*}
        \gamma_f(t) := \begin{bmatrix}
            t \\
            P_{f, i}(t)\\
            0
        \end{bmatrix}.
    \end{align*}

    Let $\kappa_f$ be the curvature (see Definition~\ref{def:curvature})
    \begin{align*}
        \kappa_f = \frac{\|\ddot{\gamma_f} \times \dot{\gamma_f}\|}{\|\dot{\gamma_f}\|^3}.
    \end{align*}

    Then, the average curvature can be expressed as
    \begin{align*}
        A_f(t_{f(i - 1)}, t_{fi}) = \frac{1}{t_{fi} - t_{f(i - 1)}} \int_{t_{f(i - 1)}}^{t_{fi}} \kappa_f(t) \, \d t.
    \end{align*}
\end{lemma}
\begin{proof}
    Although $\dom(P_{f, i})$ contains only dynamical fractal $d_{f}$ at time $t$, where $t \in (t_{f0}, t_{f1})$, we can use the function 
    \begin{align*}
        I(d_{f, t}) = t.
    \end{align*}

    For simplicity, we assume that 
    \begin{align*}
        \dom(P_{f, i}) = (t_{f0}, t_{f1}).
    \end{align*}

    Therefore, the exponential function $P_{f, i}$ can be expressed as
    \begin{align*}
        P_{f, i}(t) = a e^t,
    \end{align*}
    where $a \in \mathbb{R}$.

    Note that $P_{f, i}$ has one variable, so we let $\gamma_f(t)$ be as defined in this lemma statement. 

    By Fact~\ref{fac:curvature}, $\kappa_f$ can be expressed as shown in the lemma statement.

    Therefore, we can use $A_f(t_{f(i - 1)}, t_{fi})$ to denote its average curvature.
\end{proof}

\begin{definition}\label{def:similar_shape}
    Let everything be defined as in Lemma~\ref{lem:avergae_curvature}.

    Then, we say the two sets of exponential functions $P_f^*$ and $P_g^*$ have a similar shape if for all $i \in \{1, 2, \dots, n - 1\}$,
    \begin{align*}
        \frac{A_f(t_{f(i - 1)}, t_{fi})}{A_f(t_{fi}, t_{f(i + 1)})} \approx \frac{A_g(t_{g(i - 1)}, t_{gi})}{A_g(t_{gi}, t_{g(i + 1)})}.
    \end{align*}

    For simplicity, we let  
    \begin{align*}
        \alpha_{f, i} : = \frac{A_f(t_{f(i - 1)}, t_{fi})}{A_f(t_{fi}, t_{f(i + 1)})}.
    \end{align*}
\end{definition}

\section{Case Study: Fractal Dimension of Boston from 2000 to 2020}
\label{sec:case_study}

We start presenting our case study in this section. In Section~\ref{section:3.0}, we introduce the methods of pre-processing the satellite images. In Section~\ref{section:3.1}, we introduce the fractal dimension approximation technique: box-counting analysis. In Section~\ref{section:3.2}, we introduce another fractal dimension approximation technique: radial analysis. In Section~\ref{sub:case_study:fractal_dimension_boston}, we present our data of the fractal dimension of Boston from 2000 to 2020 approximated through the box-counting and radial analysis. In Section~\ref{sub:case_study:fractal_dimension_boroughs}, we present the data of the fractal dimension of Boston boroughs. In Section~\ref{sub:case_study:planimetric}, we introduce a planimetric map to compute the fractal dimension. In Section~\ref{sub:case_study:pros_cons}, we discuss the advantages and disadvantages of the box-counting method, radial method, satellite image, and planimetric map.

\subsection{Pre-processing the Satellite Images}
\label{section:3.0}

In this section, we introduce the techniques for pre-processing the satellite images.

We have satellite images of Boston from Google Earth. We first need to pre-process these images, either transforming them into binary images \cite{1} or using deep learning method \cite{9} to mark the human-related properties. 

\paragraph{Binary images}

The pixels in satellite images have different colors. By using Python, we transform the pixels which represent the buildings and human-related properties to white and the natural landscape to black. Thus, the fractal dimension can be computed based on the complexity of the spread of the white pixels. 

\begin{figure}[!ht]
     \centering
     \begin{subfigure}{0.49\textwidth}
         \centering
         \includegraphics[width=0.8\linewidth]{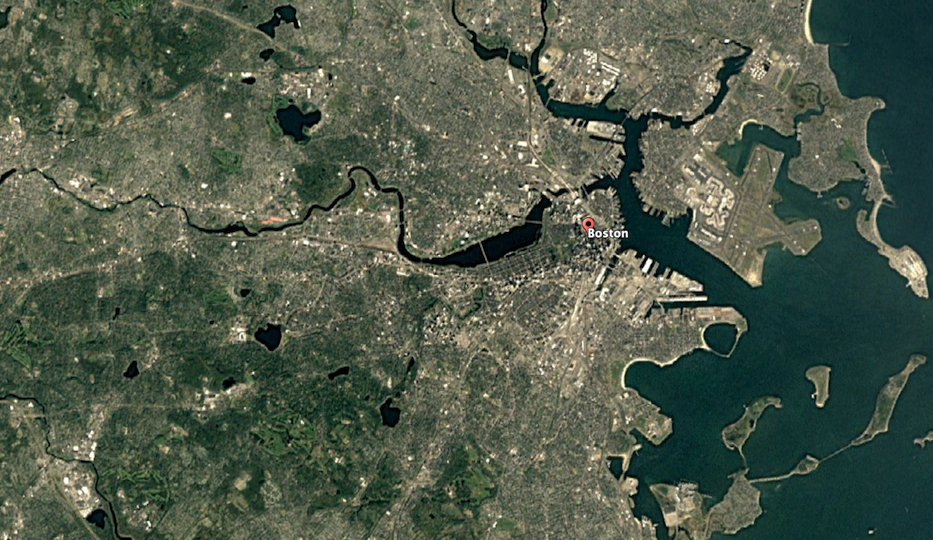}
         \caption{Satellite image of Boston region in 2000}
         \label{3}
     \end{subfigure}
     \hfill
     \begin{subfigure}{0.49\textwidth}
         \centering
         \includegraphics[width=0.8\linewidth]{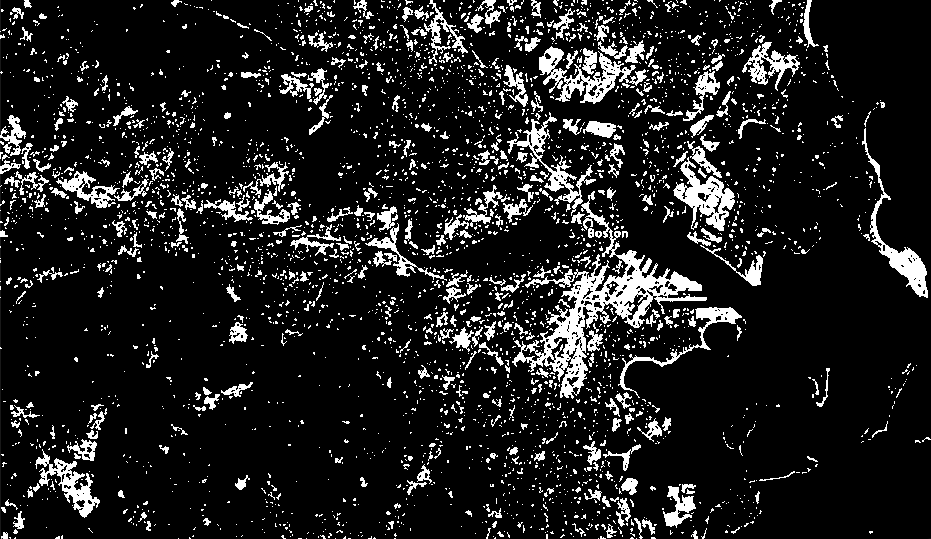}
         \caption{Binary image of Boston in 2000}
     \end{subfigure}
     \caption{Satellite images and its corresponding binary image}
\end{figure}

\paragraph{Deep learning}

The other way to pre-process is to use the deep learning method. ENVI is a software, which can learn from the color that people assign to certain pixels of satellite images, and based on that, assign the colors to similar pixels in that image. 

\begin{figure}[!ht]
     \centering
     \begin{subfigure}{0.49\textwidth}
         \centering
         \includegraphics[width=0.8\linewidth]{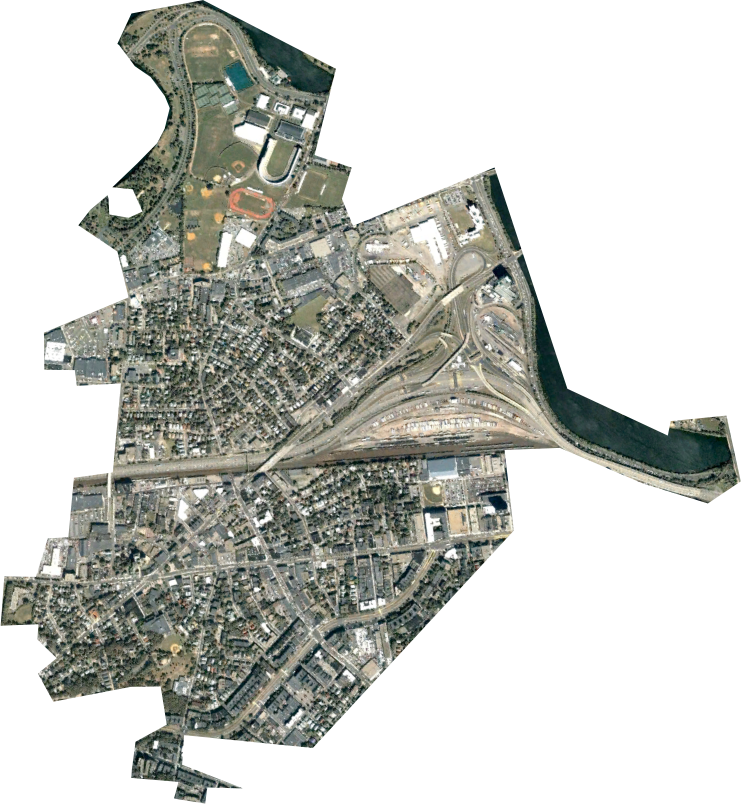}
         \caption{Satellite Image of Allston in 2002}
     \end{subfigure}
     \hfill
     \begin{subfigure}{0.49\textwidth}
         \centering
         \includegraphics[width=0.8\linewidth]{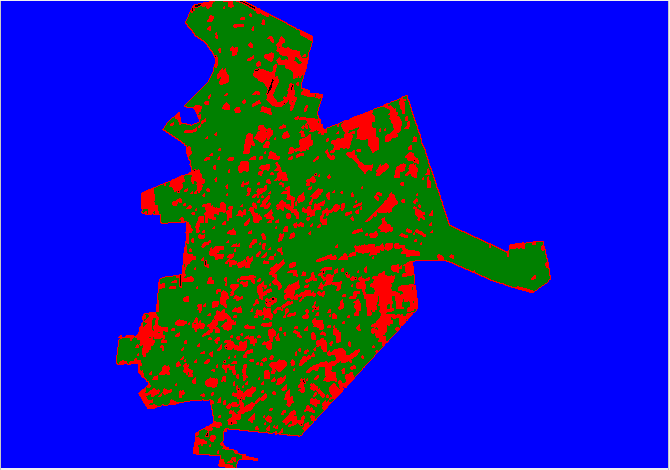}
         \caption{Satellite Image of Allston in 2002 processed by ENVI}
     \end{subfigure}
        \caption{Image of Allston in 2002 before and after processing by ENVI}
        \label{9}
\end{figure}

In Figure~\ref{9}, we assigned the red color to a certain number of buildings and yards, the green color to several vegetation and roads, and the blue color to the background. Since the buildings and yards are white and light gray, vegetation is green, and roads are dark gray, ENVI can learn from that and transform all the remaining buildings, roads, and yards to the corresponding colors. 

\subsection{Box-counting Analysis}
\label{section:3.1}

We first cut the specific area that we would like to analyze.

\begin{figure}[!ht]
    \centering
    \includegraphics[width=0.45\linewidth]{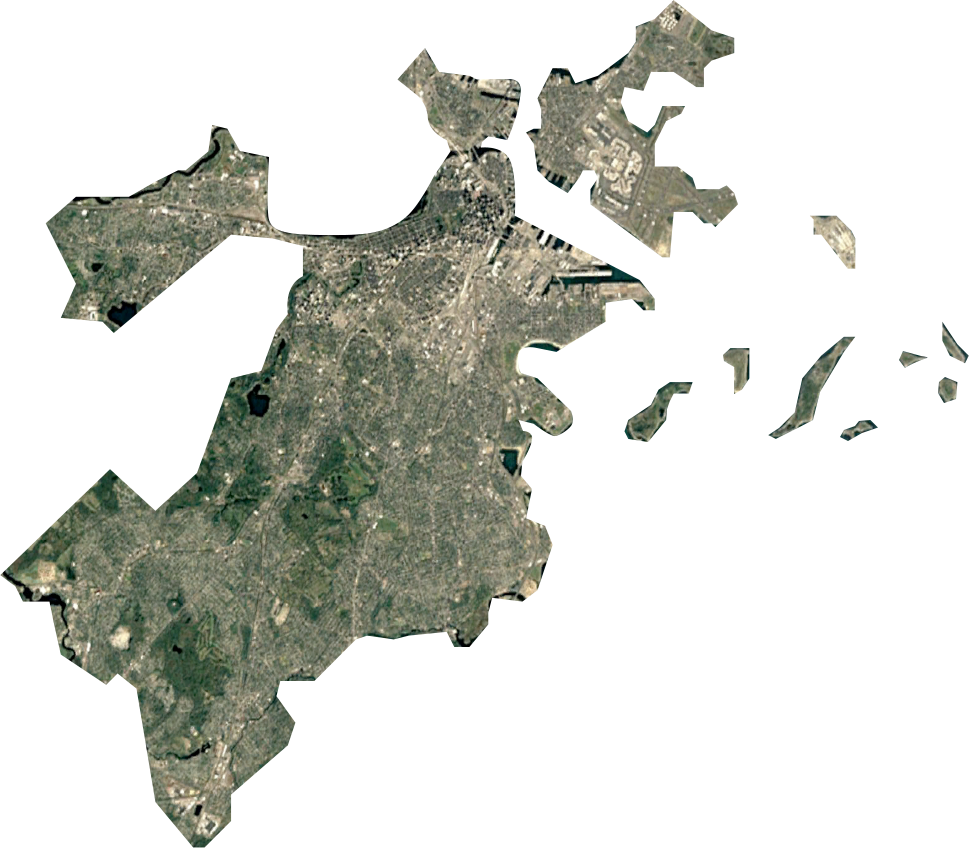}
    \caption{Boston satellite image in 2000}
    \label{5}
\end{figure}

\begin{definition}[Box-counting analysis \cite{7}]
    The box-counting analysis is a grid analysis that the rectangular girds covering a certain image with the side length $\epsilon$, and for all $\epsilon > 0$, there exists a $N(\epsilon)$, such that:
    \begin{align*}
        D = \lim_{\epsilon \to 0} \frac{\log{(N(\epsilon))}}{-\log{(\epsilon)}}.
    \end{align*}
\end{definition}

$N(\epsilon)$, or count, represents the total number of the grids which contain a certain part of the object. 

\begin{figure}[!ht]
    \centering
    \begin{subfigure}[b]{0.49\textwidth}
        \centering
        \includegraphics[width=0.8\linewidth]{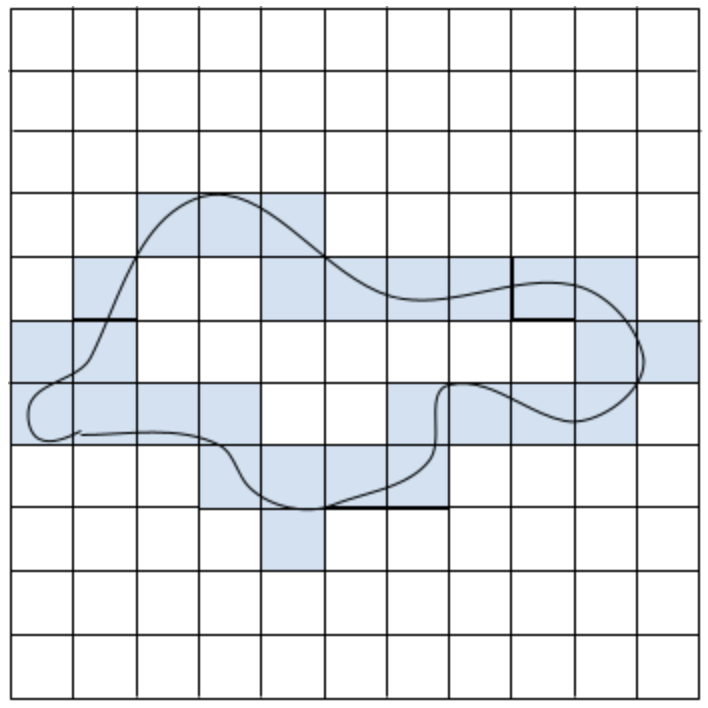}
        \caption{Boxes with higher side lengths}
        \label{bcd}
    \end{subfigure}
    \hfill
    \begin{subfigure}[b]{0.49\textwidth}
        \centering
        \includegraphics[width=0.8\linewidth]{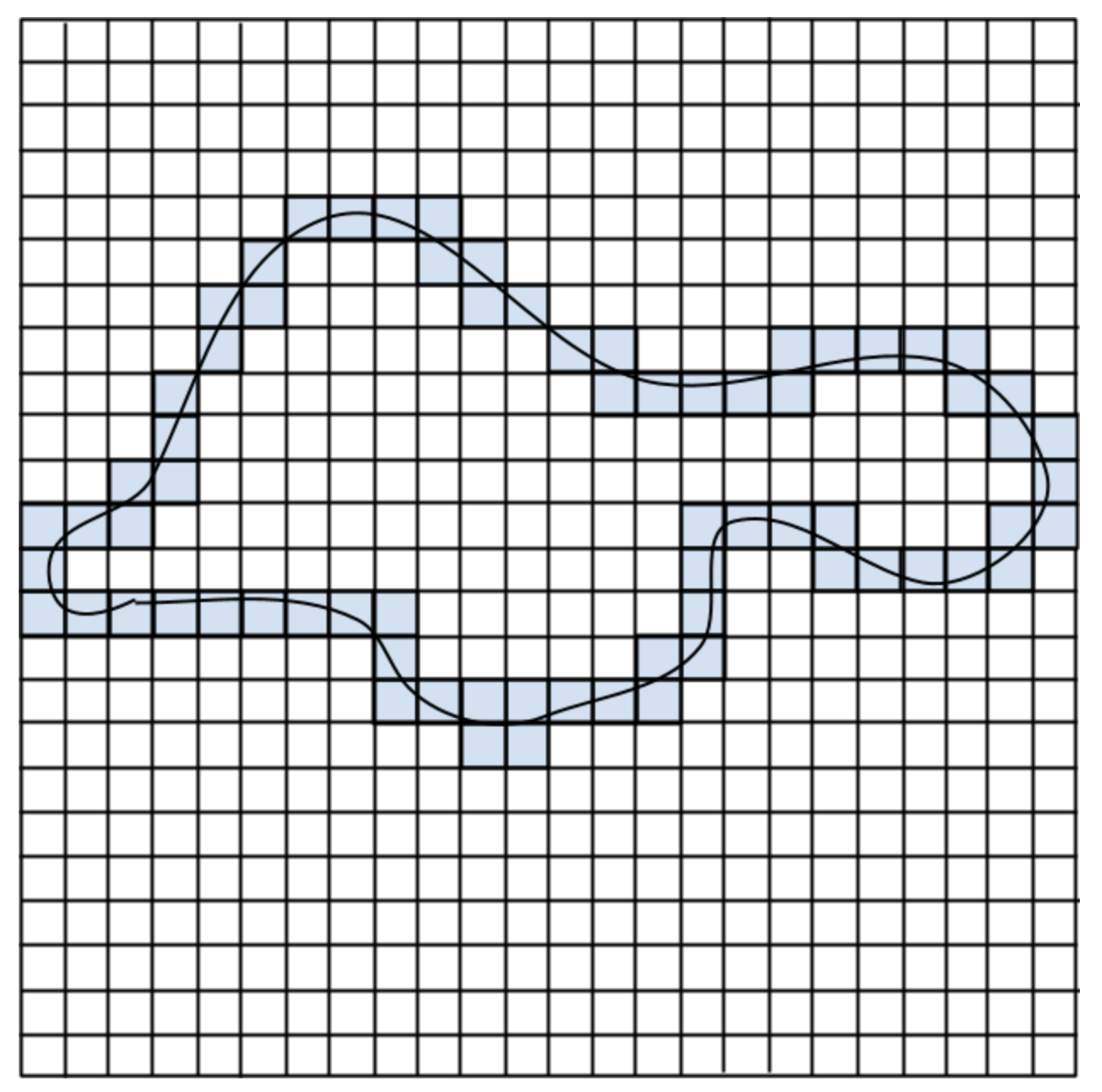}
        \caption{Boxes with lower side lengths}
        \label{bcx}
    \end{subfigure}
    \caption{Number of grids and side length}
\end{figure}

\begin{figure}[!ht]
    \centering
    \begin{subfigure}[b]{0.49\textwidth}
        \centering
        \includegraphics[width=0.8\linewidth]{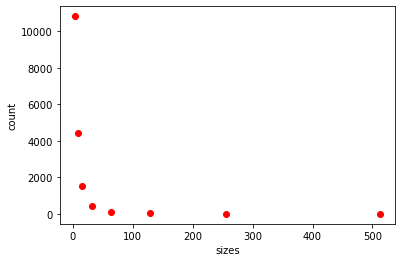}
        \caption{Different grid sizes and their corresponding counts}
        \label{5a}
    \end{subfigure}
    \hfill
    \begin{subfigure}[b]{0.49\textwidth}
        \centering
        \includegraphics[width=0.8\linewidth]{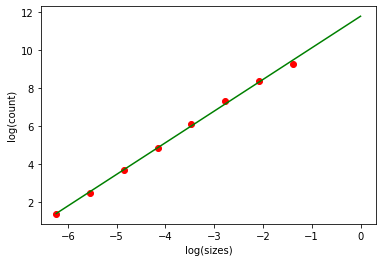}
        \caption{The linear regression line of (a)}
        \label{5b}
    \end{subfigure}
        \caption{Fractal Dimension of Boston in 2000 approximated by Box-counting Method}
\end{figure}

For example, in Figure~\ref{bcd} and Figure~\ref{bcx}, the number of colored grids $N(\epsilon)$ increases as the side length $\epsilon$ decreases. We keep changing $\epsilon$ until we get a pattern similar to Figure~\ref{5a}. 

When starting to analyze the satellite image of Boston in 2000 by using the box-counting method, we get the following result. Figure \ref{5a} can be modeled by Eq.~\eqref{eq1}, and Figure \ref{5b} is modeled by Eq.~\eqref{eq:2.2}, which is the linearized equation of Eq.~\eqref{eq1}. The slope of the linear equation in Figure \ref{5b} is the fractal dimension approximated by the box-counting method, which is equal to $1.475$. 

\subsection{Radial Analysis}
\label{section:3.2}

Radial analysis, on the other hand, refers to a specific point known as a counting center. It measures the law of distribution around that point \cite{14}. First, we find the center of the image and draw a small circle around that center. Let $r$ be the radius of that circle. We change $r$. For each $r$, the total number of the occupied point $N$ is counted. Their relationship can be expressed by:
\begin{align*}
        N &= R^D\\
        \log(N) &= D\log(R).
\end{align*}

Similarly, the function in Figure~\ref{8b} can also be linearized, and its slope, $D$, is the fractal dimension approximated by the radial method. 

\begin{figure}[!ht]
     \centering
     \begin{subfigure}{0.49\textwidth}
         \centering
         \includegraphics[width=0.8\linewidth]{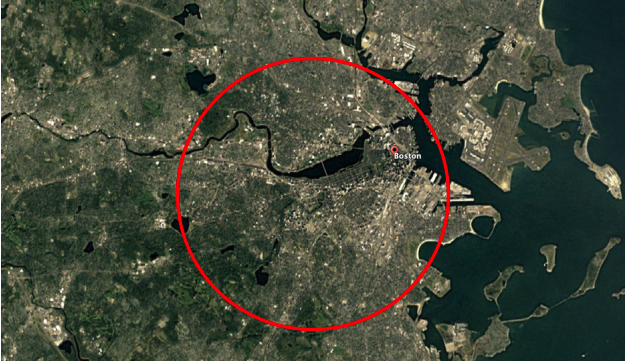}
         \caption{Radial method with the uncut Boston satellite image}
     \end{subfigure}
     \hfill
     \begin{subfigure}{0.49\textwidth}
         \centering
         \includegraphics[width=0.7\linewidth]{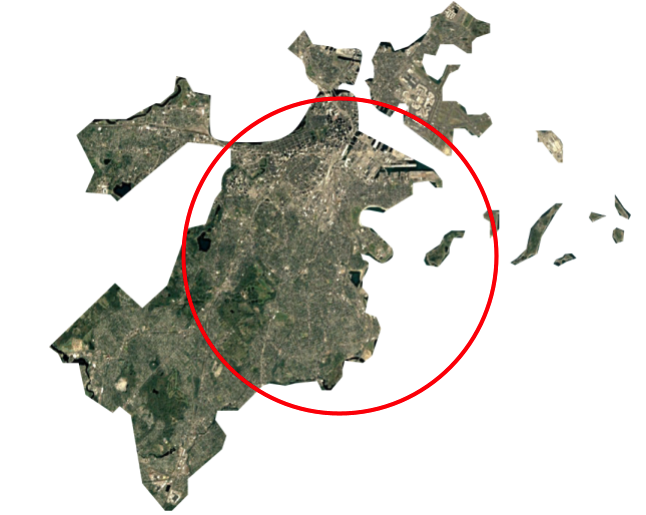}
         \caption{Radial method with the cut Boston satellite image}
         \label{7b}
     \end{subfigure}
     \caption{Satellite images in radial method}
\end{figure}

However, we cannot use the cut satellite images while using the radial method. From Figure~\ref{7b}, if we cut out the Boston area from the satellite image, and when the radius of the circle is large, it will contain the background region (the white area). Since there are no buildings in those regions, the relationship between the occupied point $N$ and the radius $R$ will be inaccurate. 

\begin{figure}[!ht]
     \centering
     \begin{subfigure}{0.49\textwidth}
         \centering
         \includegraphics[width=0.8\linewidth]{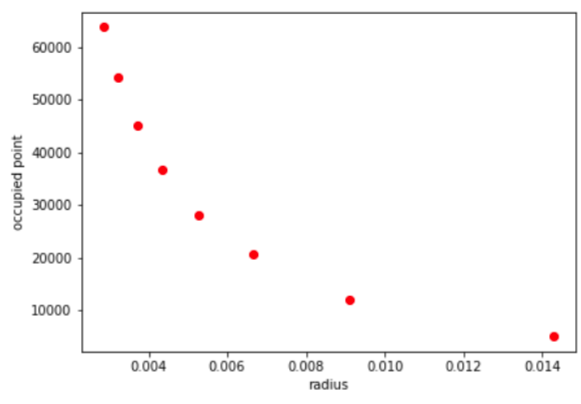}
         \caption{Different radius and their corresponding occupied point}
         \label{8a}
     \end{subfigure}
     \hfill
     \begin{subfigure}{0.49\textwidth}
         \centering
         \includegraphics[width=0.8\linewidth]{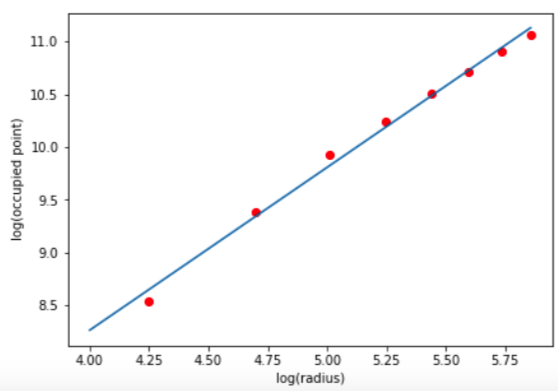}
         \caption{The linear regression line of (a)}
         \label{8b}
     \end{subfigure}
        \caption{Fractal Dimension of Boston in 2000 approximated by radial method}
\end{figure}

By looking for the slope of the function in Figure~\ref{8b}, we get that the fractal dimension is equal to $1.541$, approximated by the radial method.

\subsection{Fractal Dimension of Boston Approximated by Box-counting Method and Radial Method}
\label{sub:case_study:fractal_dimension_boston}

\begin{figure}[!ht]
    \centering
    \includegraphics[width=0.45\linewidth]{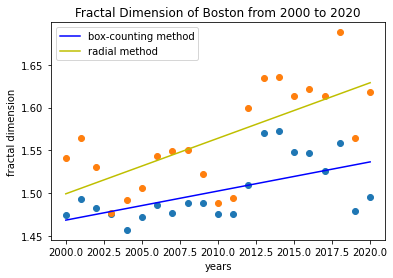}
    \caption{Fractal dimension of Boston from 2000 to 2020}
    \label{fig:FDBR}
\end{figure}

We use the box-counting method in Section \ref{section:3.1} and the radial method in Section \ref{section:3.2} to analyze the binary images of Boston from 2000 to 2020. By using Python, we get Figure~\ref{fig:FDBR}.

We can get that the fractal dimension approximated by using the box-counting method is generally lower than that generated by the radial method, and it has a lower standard deviation\footnote{The standard deviation of the fractal dimension of Boston from 2000-2020 approximated by the box-counting method is 0.0348, and that of the radial method is 0.0573.}.

\subsection{Fractal Dimension of Boston's Boroughs}
\label{sub:case_study:fractal_dimension_boroughs}

In this section, we present the fractal dimension of Boston's boroughs.

It is not accurate enough to solely analyze the general satellite images of Boston because the diameter of Boston is large, and the degree of compression cause the satellite image to omit some details, leading us to underestimate the fractal dimension. There are two more ways of adding accuracy to the approximated fractal dimension. First, we analyze the fractal dimension of the Boston boroughs, the smaller administrative units than cities, like Brooklyn in New York City. By Theorem~\ref{thm:self_similar}, Boston boroughs are expected to have similar fractal dimensions as Boston.

Compared to the satellite images of Boston, the satellite images of Boston boroughs have a higher spatial resolution but are more susceptible to clouds and brightness, so we use the deep learning method (see Section~\ref{section:3.0}) to process them to eliminate clouds and brightness issue. 

We analyzed two Boston boroughs: Allston and Beacon Hill. Here is their fractal dimension:

\begin{table}[!ht]
\begin{center}
\begin{tabular}{ |c|c|c|c| } 
 \hline
 year & Allston & Beacon Hill & year \\ 
 \hline
 2000 & 1.421 & 1.697 & 2001\\ 
 \hline
 2002 & 1.525 & 1.693 & 2004\\ 
 \hline
 2013 & 1.527 & 1.686 & 2015\\ 
 \hline
 2020 & 1.542 & 1.688 & 2019\\ 
 \hline
\end{tabular}
\caption{Fractal Dimension approximated by the box-counting method of Allston and Beacon Hill}
\end{center}
\end {table}

\begin{figure}[!ht]
    \centering
    \includegraphics[width=0.45\linewidth]{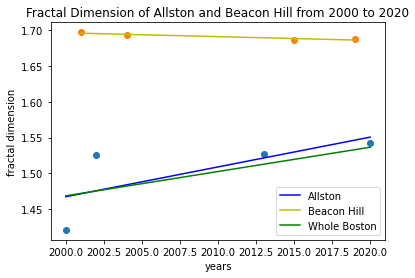}
    \caption{Fractal Dimension approximated by the box-counting method of Allston and Beacon Hill from 2000 to 2020}
    \label{fig:BHD}
\end{figure}

Beacon Hill is located in Boston downtown, so it is more urbanized. It has a higher fractal dimension than Boston since it had already been very urbanized since 2000. We use that as the upper bound of the fractal dimension of Boston because the whole city cannot be more urbanized than its downtown area. Allston, on the other hand, is located in the west part of Boston. It was continuously developed throughout these years, so it has a very similar fractal dimension to that of Boston. This result is the same as Theorem \ref{thm:3.6}.

\subsection{Approximating the Fractal Dimension Using Planimetric Images}
\label{sub:case_study:planimetric}

In this section, we approximate fractal dimensions by planimetric images.

Analyzing the fractal dimension of the planimetric images \cite{18} of Boston buildings is another way of learning its complexity. Unlike the satellite images, even though the planimetric images only include the buildings of Boston, which may not thoroughly reflect the degree of urbanization or the complexity of Boston, they are not affected by weather, brightness, seasons, and clouds. 

\begin{figure}[!ht]
    \centering
    \includegraphics[width=0.5\linewidth]{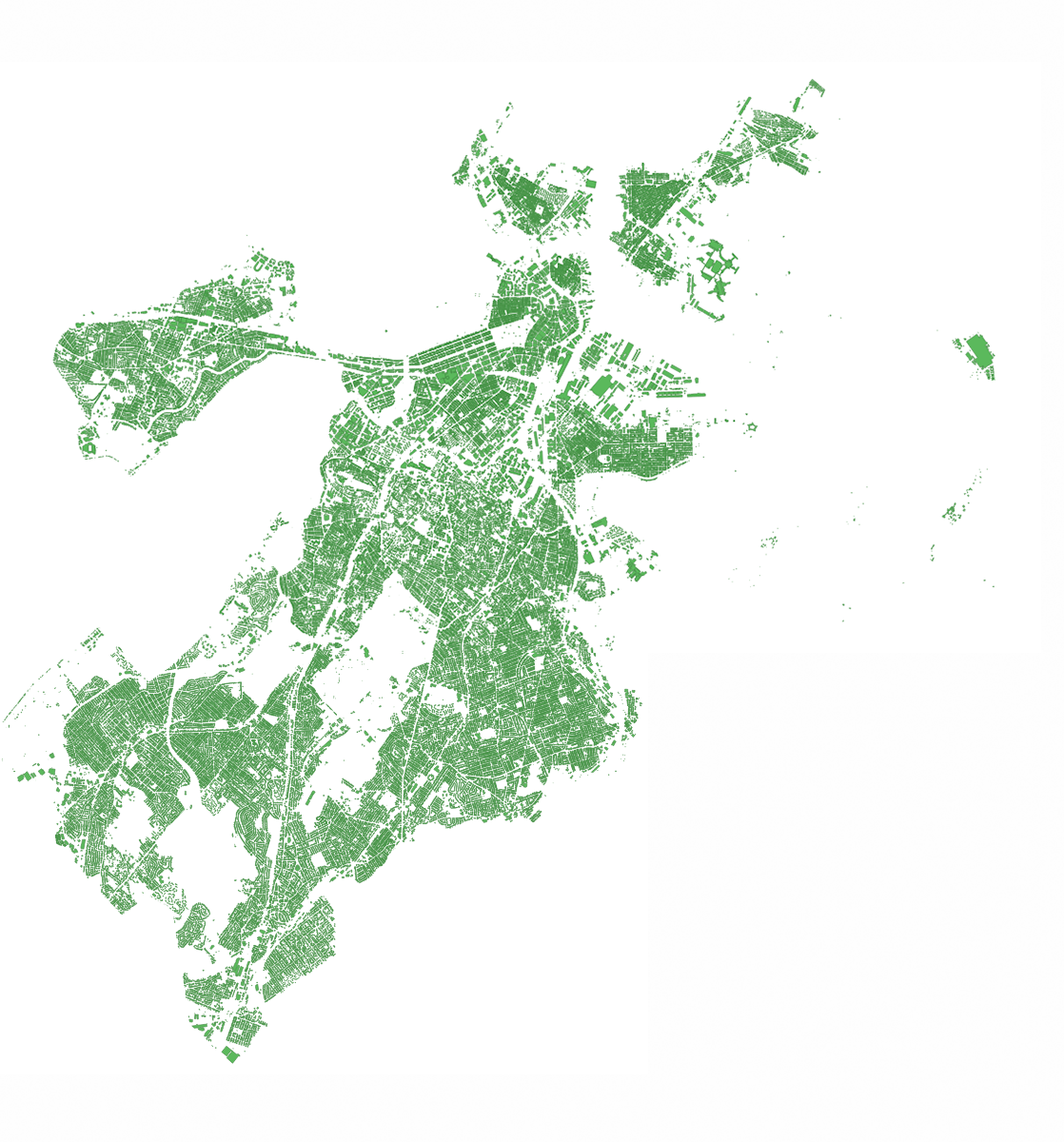}
    \caption{Planimetric map of Boston in 2020}
\end{figure}

Similar to the images processed by ENVI, the planimetric map can only be approximated by the box-counting method. By using it, we get the fractal dimension of the planimetric map of Boston in 2020 is equal to $1.517$. 

\subsection{Strengths and Weaknesses of Satellite Images, ENVI, and Planimetric Map}
\label{sub:case_study:pros_cons}

In this section, we analyze the strengths and weaknesses of satellite images, ENVI, and the planimetric map.

The advantage of satellite images is that they can include some small details of the city because they are real images taken by the satellites. However, when we transform them into binary images, the errors are generated because we have to find a certain value\footnote{The lower that value is, the more pixels will be transformed into white.}, which tells the code what kinds of colors we would like to become white. However, it is hard to find a precise value, which can transform all the properties in the image correctly by using our eyes. Also, when observing the whole city, we realize that the spatial resolution is relatively low, and when observing the boroughs of Boston, even though we find that the satellite images have higher spatial resolutions, the images are susceptible to cloud, fog, and brightness.

Being processed by ENVI, those satellite images will not have these issues because we can tell that software which kinds of pixels we would regard as human-related properties, and which kinds of pixels we would regard as nonhuman-related properties. Based on the examples that we give it, it can differentiate other pixels. Another unexpected advantage of ENVI is that since the heights of the buildings are also an indicator of how urbanized the city is, ENVI can take that as a variable when we process it. However, even if we give some examples of what we want, it cannot always color every pixel precisely, especially when it tries to identify some similar colors, like light gray (buildings) and dark gray (roads). Furthermore, it is very time-consuming to process these images because ENVI spends a very long time calculating and requires a relatively high computer setup.

Planimetric maps, on the other hand, are the most accurate way of estimating the complexity of Boston buildings because it is made by experts based on the real sizes of the Boston buildings, and it is not influenced by bad weather, errors of binary image and ENVI, brightness, and low spatial resolution. The fractal dimension generated by the planimetric map is very similar to that generated by the satellite image in 2020\footnote{By using the planimetric map, we get the fractal dimension which is equal to 1.517, and by using the satellite image, we get that the fractal dimension is equal to 1.496.}. Nonetheless, we only have the planimetric map of Boston in 2020 from \cite{18}.

\section{The Analysis of the Trend of the Fractal Dimension of Boston}
\label{sc:atfdb}

In this section, we apply the properties that we developed in Section~\ref{sc:df} to our case study.

In Section~\ref{sub:atfdb:difference}, we model the short-term data by the difference equation. In Section~\ref{sc:5.2}, we compare the population growth of Boston and Manhattan. In Section~\ref{sub:atfdb:differential}, we model the long-term data by the differential equation.

\subsection{Fitting the Data with the Difference Equation}
\label{sub:atfdb:difference}

In this section, we use a difference equation to extrapolate the short-term trend of the fractal dimension. 

Based on Definition~\ref{def:difference}, we let 
\begin{align}\label{eq:difference}
    \Phi^*(d_{f, t_m}) = c_1 t_m + c_2 + \frac{(t_m - c_3)^2 \sin(t_m - c_4)}{c_5},
\end{align}
where for $i \in \{1, 2, 3, 4, 5\}$, $c_i \in \mathbb{R}$ and $t_m$ represents time. 

In the short run, we predict that the growth of the fractal dimension is linear because regardless of what the long-term trends will be, the local trend can be regarded as a linear function, so we use $c_1 t_m + c_2$ to represent it and $\frac{(t_m - c_3)^2 \sin(t_m - c_4)}{c_5}$ represents the oscillation of this trend. In order to get $c_1$, $c_2$, $c_3$, $c_4$, and $c_5$, which makes the prediction closer to the actual fractal dimension, we calculate the error and make it as small as possible, namely
\begin{align*}
    \min_{\forall i, c_i \in \mathbb{R}} (\sum_{t_m = 2000}^{2020}  |\Phi^*(d_{f, t_m}) - \Phi(d_{f, t_m})|),
\end{align*}
where $\Phi^*$ is Eq.~\eqref{eq:difference} and $\Phi$ is the fractal dimension approximated by the box-counting analysis (see Figure~\ref{fig:FDBR}).

\begin{figure}[!ht]
     \centering
     \begin{subfigure}{0.49\textwidth}
         \centering
         \includegraphics[width=0.8\linewidth]{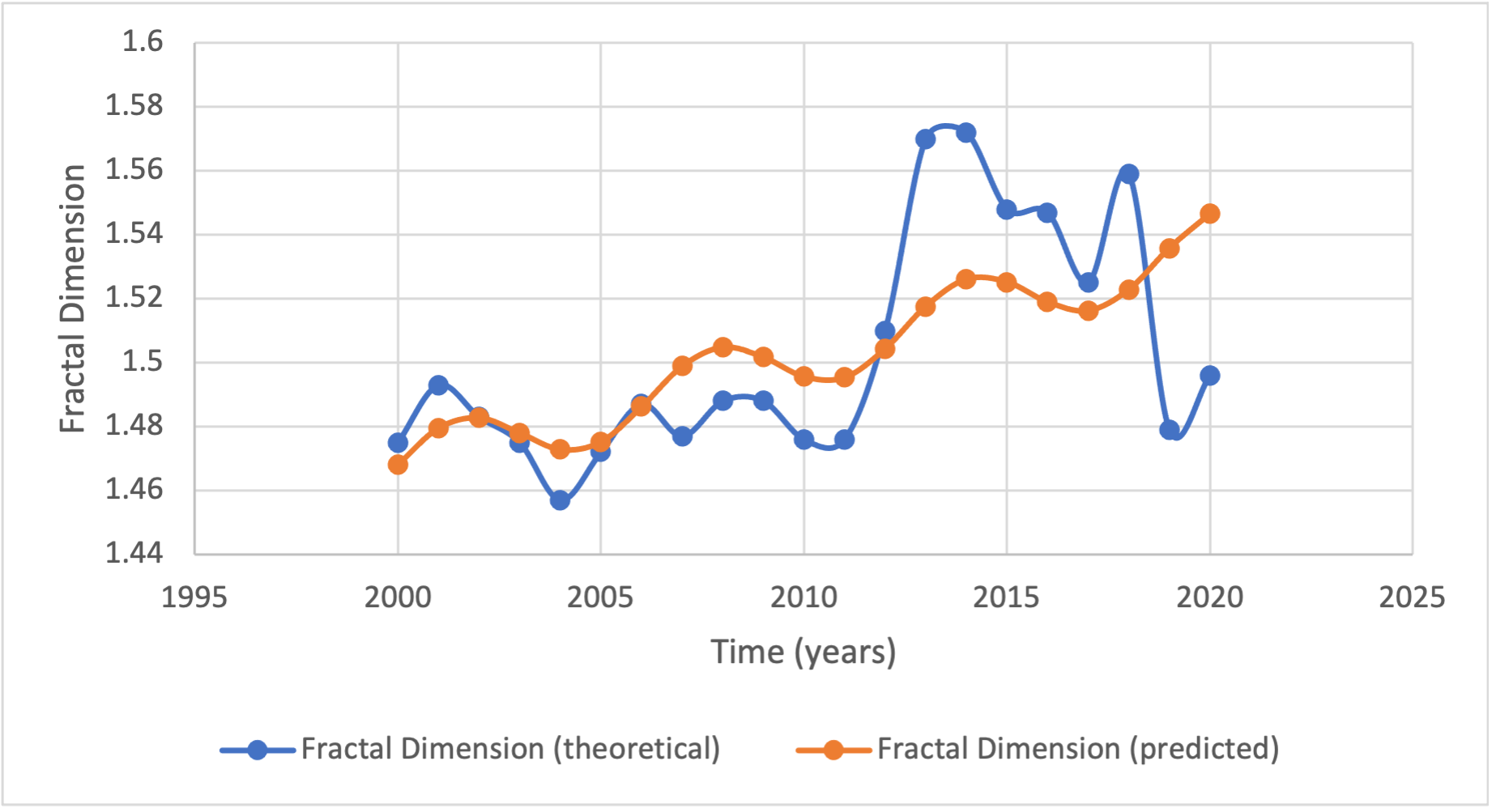}
     \end{subfigure}
     \hfill
     \begin{subfigure}{0.49\textwidth}
        \centering
        \includegraphics[width=0.8\linewidth]{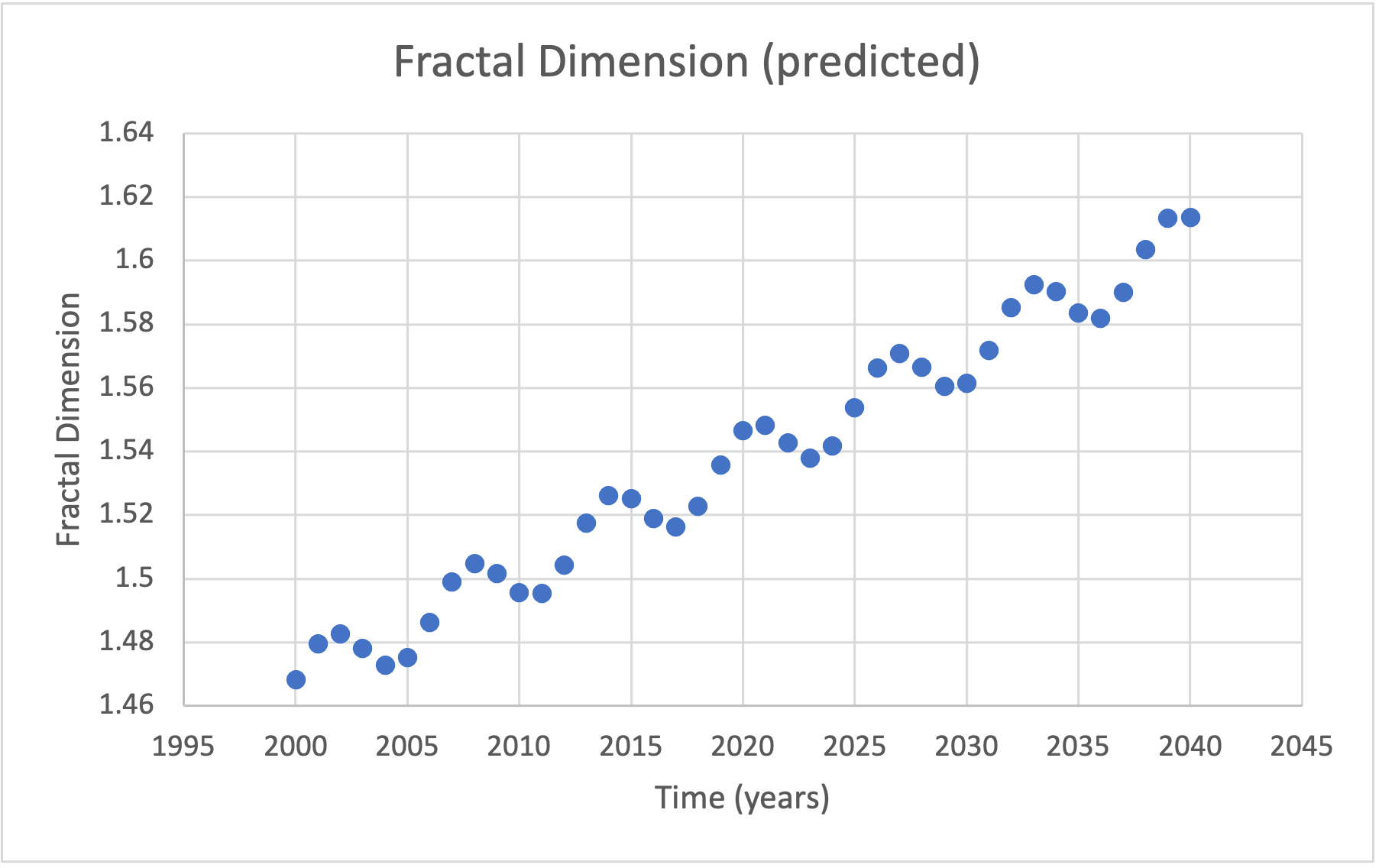}
     \end{subfigure}
        \caption{Modeling with difference equation}
\end{figure}

After getting $c_1$ through $c_5$, we plug them into the equation:
\begin{align}\label{eq:difference:c}
    \Phi^*(d_{f, t_m}) = 0.003481t_m + 5.494 + \frac{(t_m - 659.6)^2 \sin(t_m - 1.847)}{180000000}.
\end{align}

Then, we use some simple algebra to transform Eq.~\eqref{eq:difference:c} into the form of Definition~\ref{def:difference}:
\begin{align*}
    & ~ \Phi^*(d_{f, t_{m + 1}}) - \Phi^*(d_{f, t_{m}}) \\
    = & ~ 0.003481+\frac{(t_m - 658.6)^2\sin(t_m - 0.847) - (t_m - 659.6)^2 \sin(t_m - 1.847)}{180000000}.
\end{align*}

\subsection{Comparing the Population growth of Boston and Manhattan}
\label{sc:5.2}

In this section, we compare the population growth of Boston and Manhattan by analyzing the curvature of the exponential functions.

We only get the fractal dimension of Boston from 2000-2020. To extrapolate the long-term trend, we want to show that the long-term trend of Manhattan is similar to that of Boston (see Definition~\ref{def:population}). 

\begin{figure}[!ht]
     \centering
     \begin{subfigure}{0.49\textwidth}
         \centering
         \includegraphics[width=0.8\linewidth]{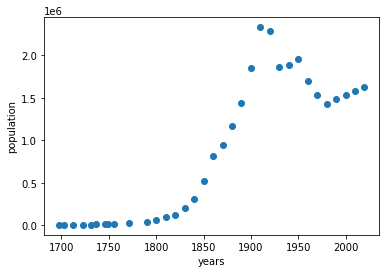}
         \caption{Population of Manhattan (the data is from \cite{27})}
         \label{14a}
     \end{subfigure}
     \hfill
     \begin{subfigure}{0.49\textwidth}
        \centering
        \includegraphics[width=0.8\linewidth]{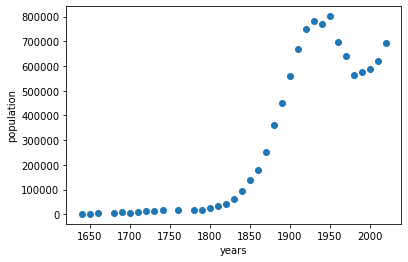}
        \caption{Population of Boston (the data is from \cite{28})}
        \label{14b}
     \end{subfigure}
        \caption{Comparison between the population in Boston and Manhattan}
        
        \label{14}
\end{figure}

\pagebreak

\begin{figure}[!ht]
     \centering
     \begin{subfigure}{0.44\textwidth}
         \centering
         \includegraphics[width=1\linewidth]{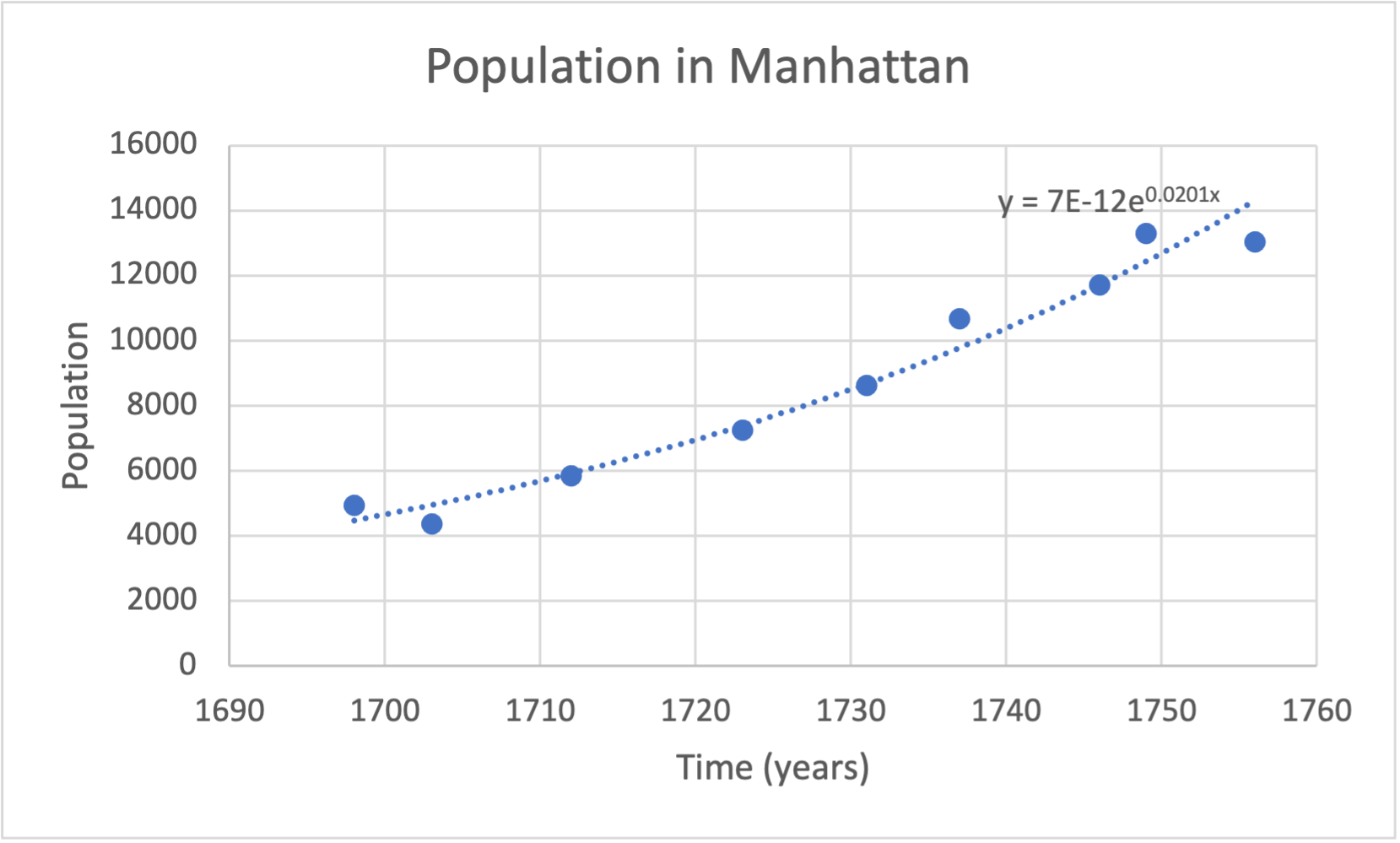}
         \includegraphics[width=1\linewidth]{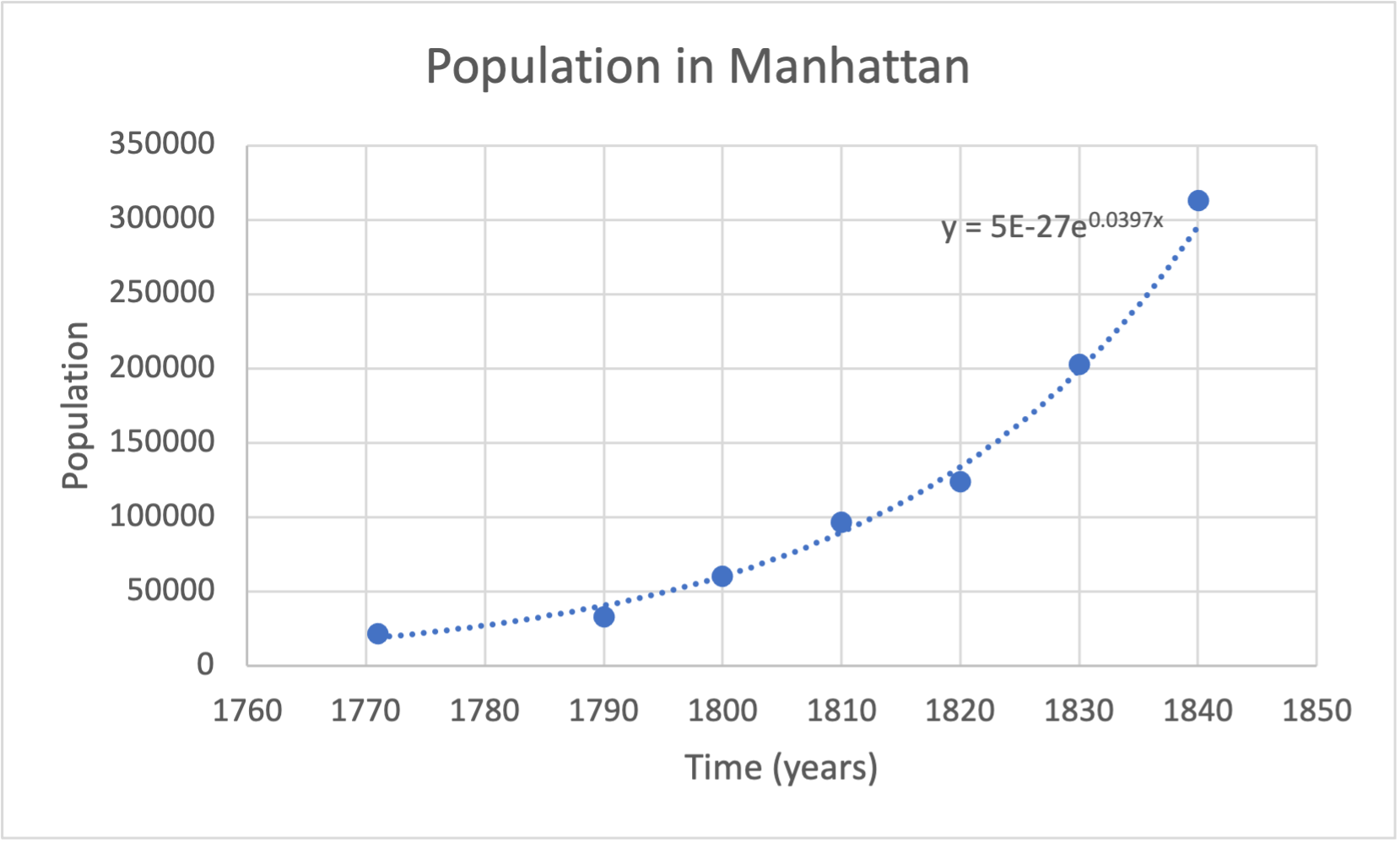}
         \includegraphics[width=1\linewidth]{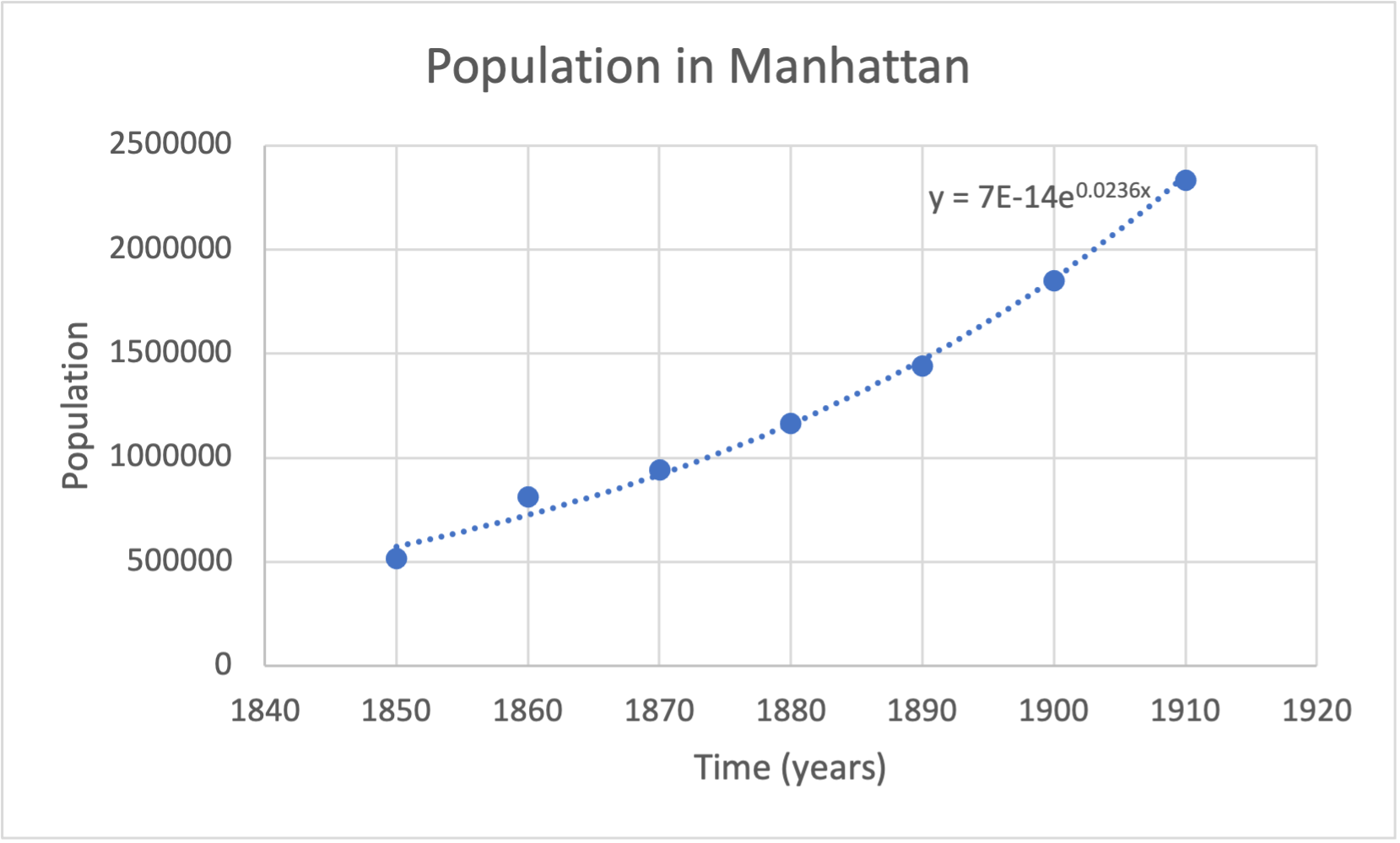}
         \includegraphics[width=1\linewidth]{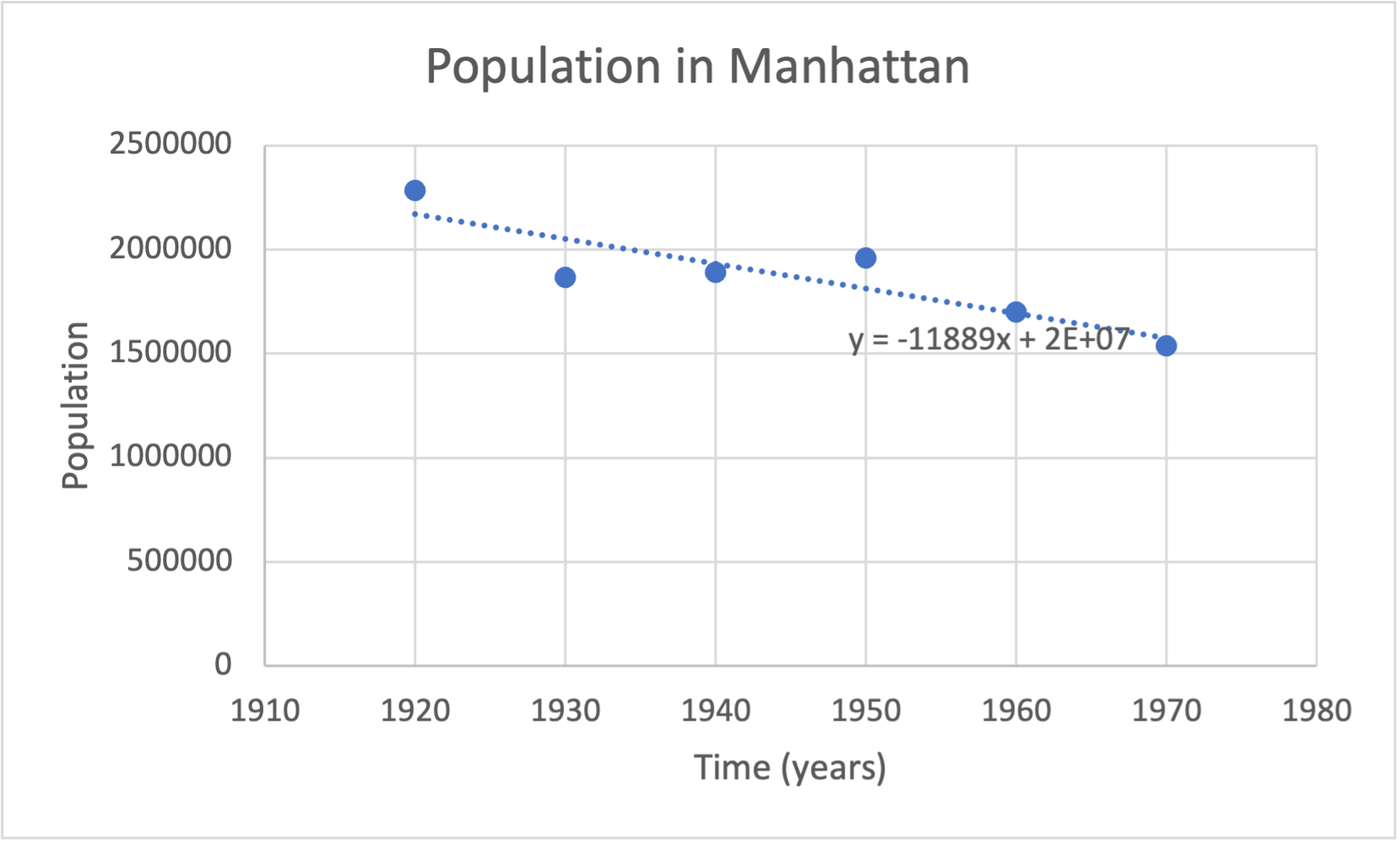}
         \includegraphics[width=1\linewidth]{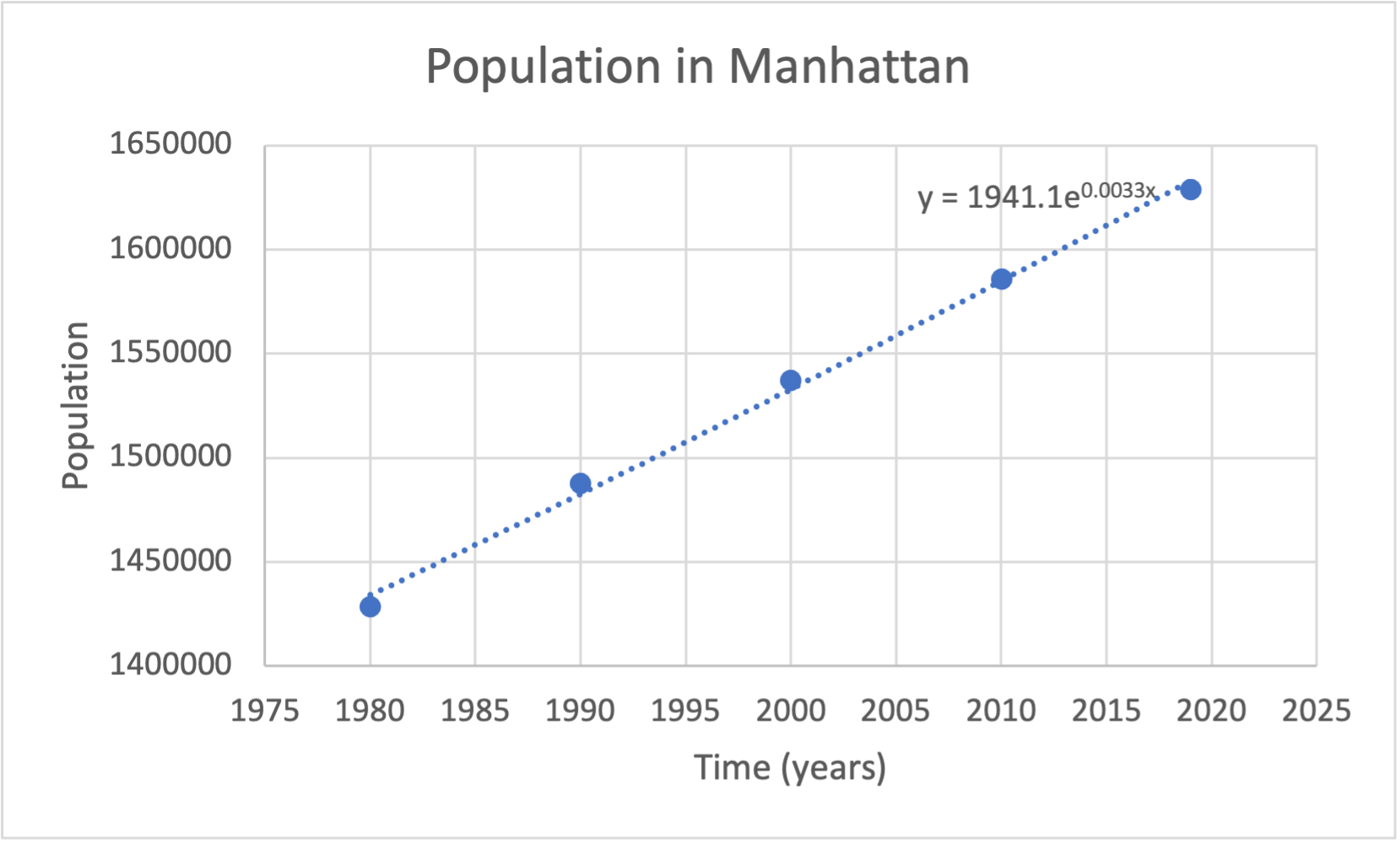}
         \caption{Population of Manhattan}
     \end{subfigure}
     \hfill
     \begin{subfigure}{0.44\textwidth}
        \centering
        \includegraphics[width=0.915\linewidth]{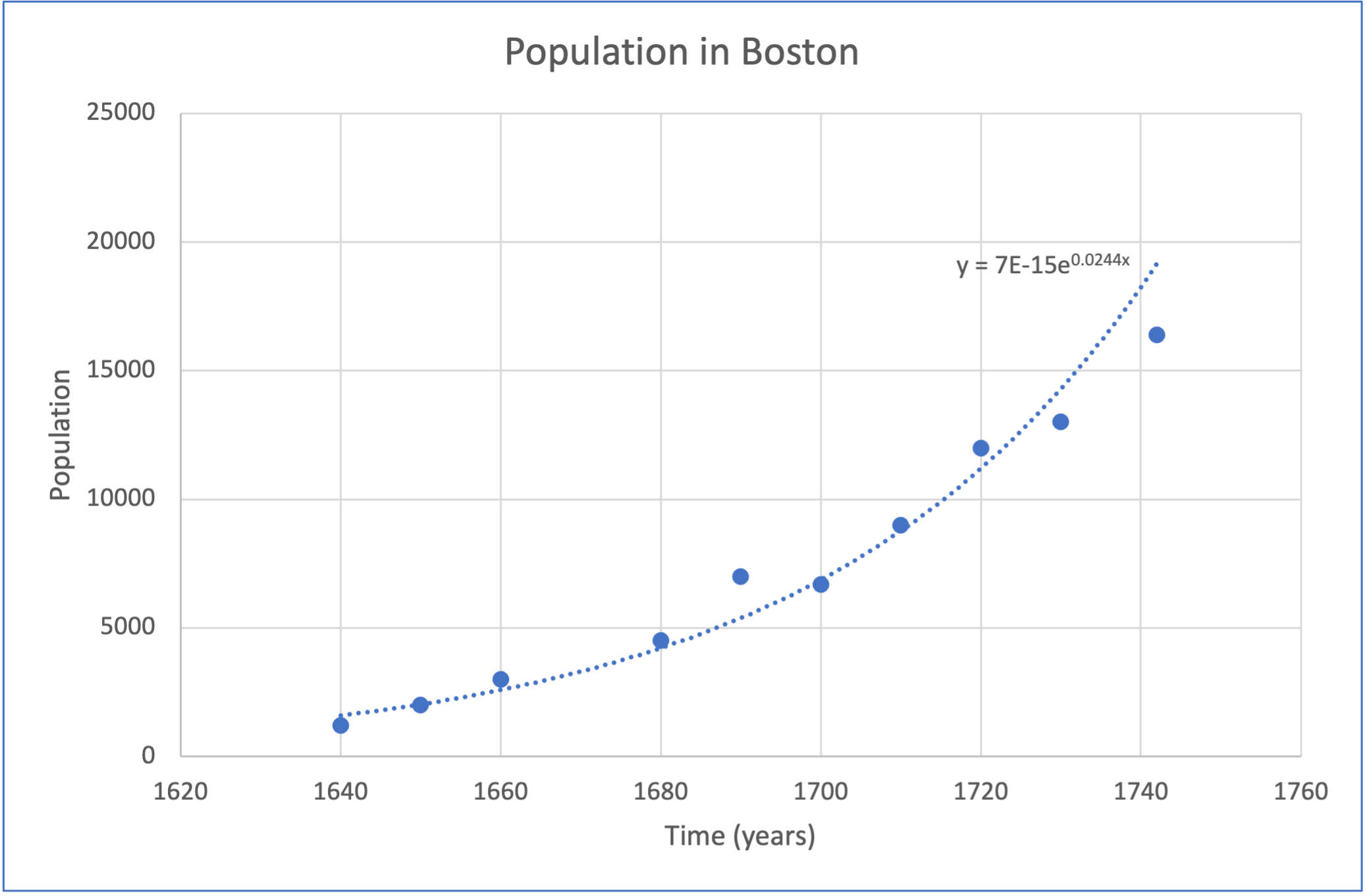}
        \includegraphics[width=0.915\linewidth]{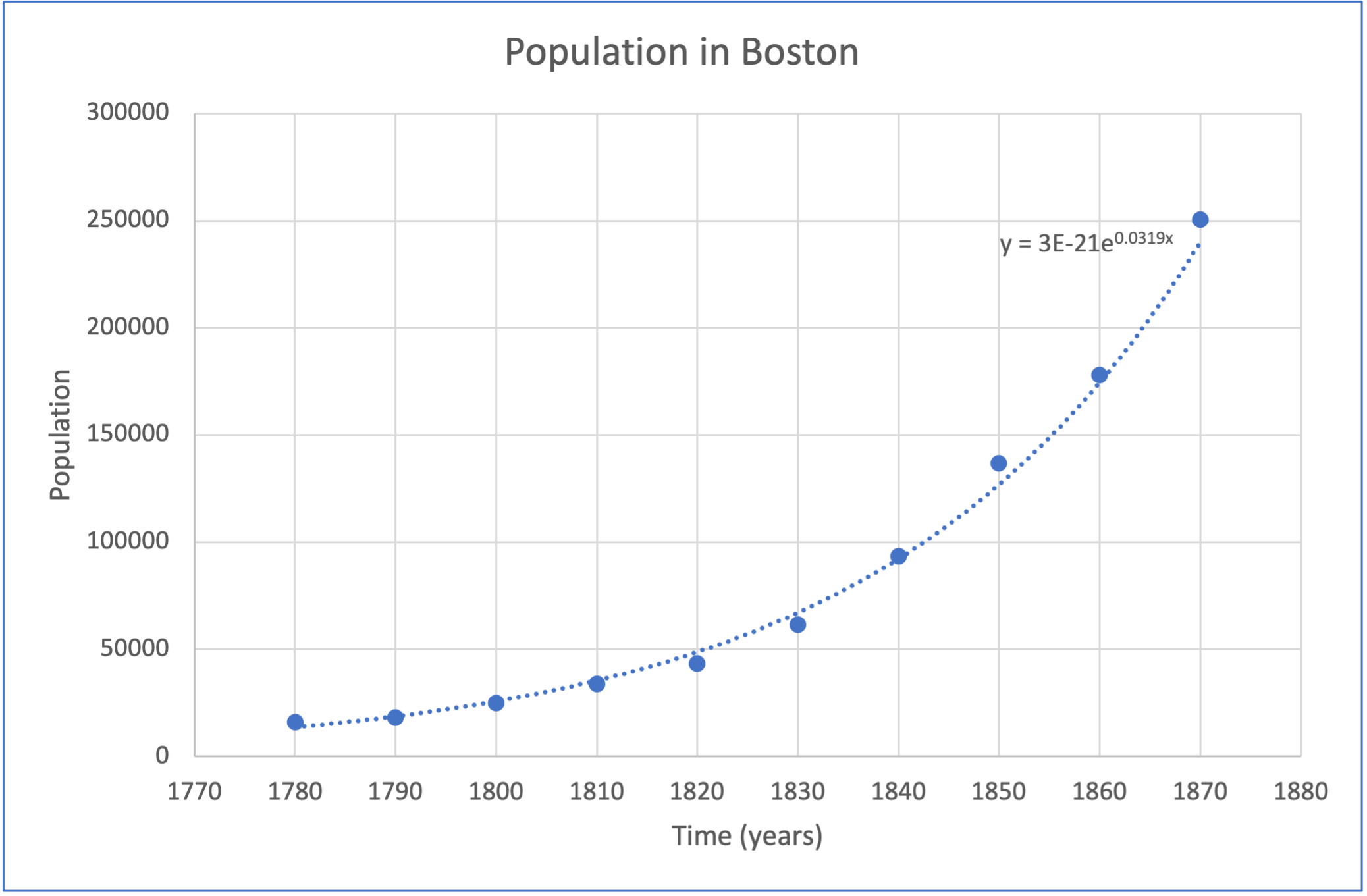}
        \includegraphics[width=0.915\linewidth]{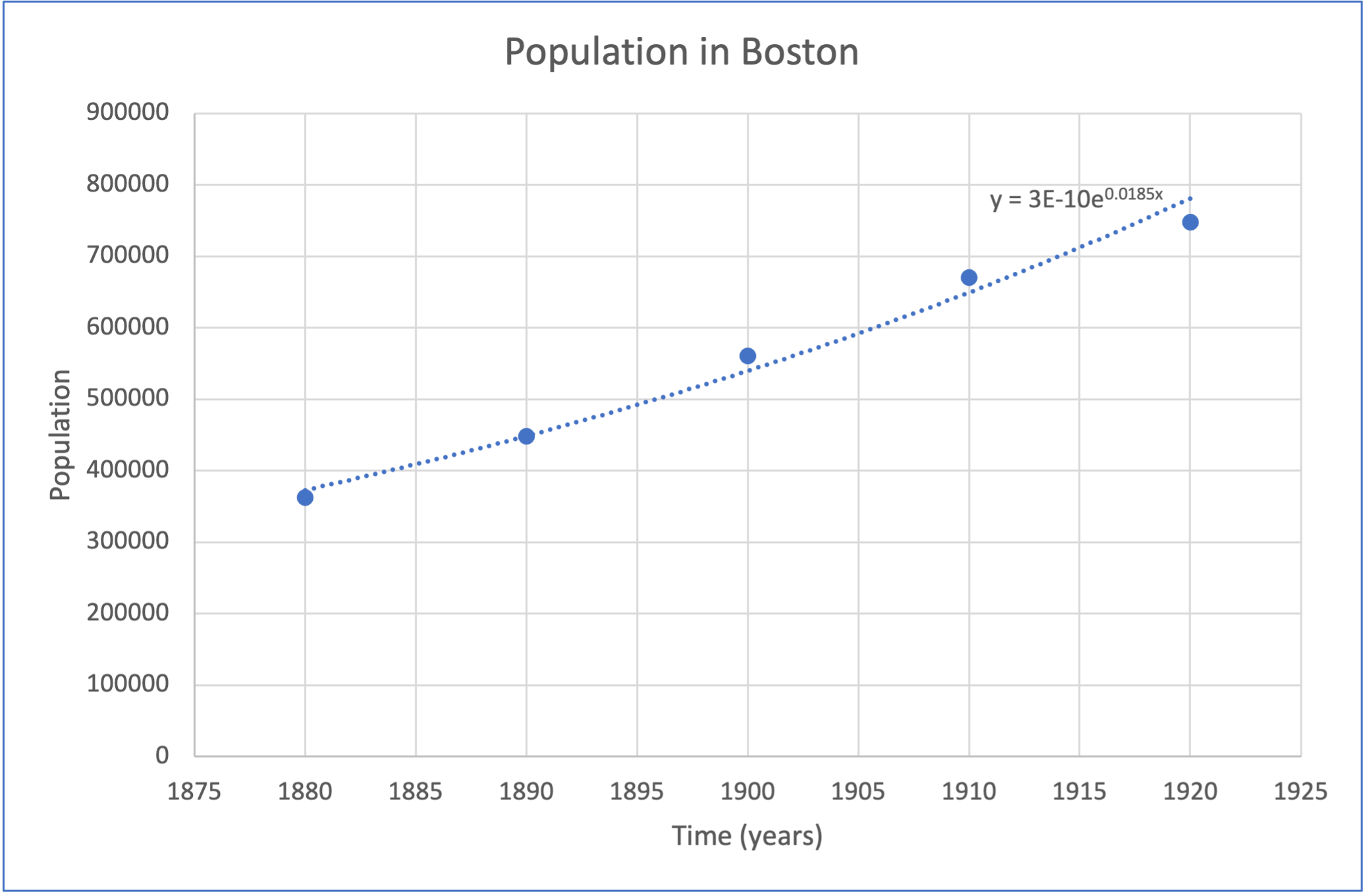}
        \includegraphics[width=0.915\linewidth]{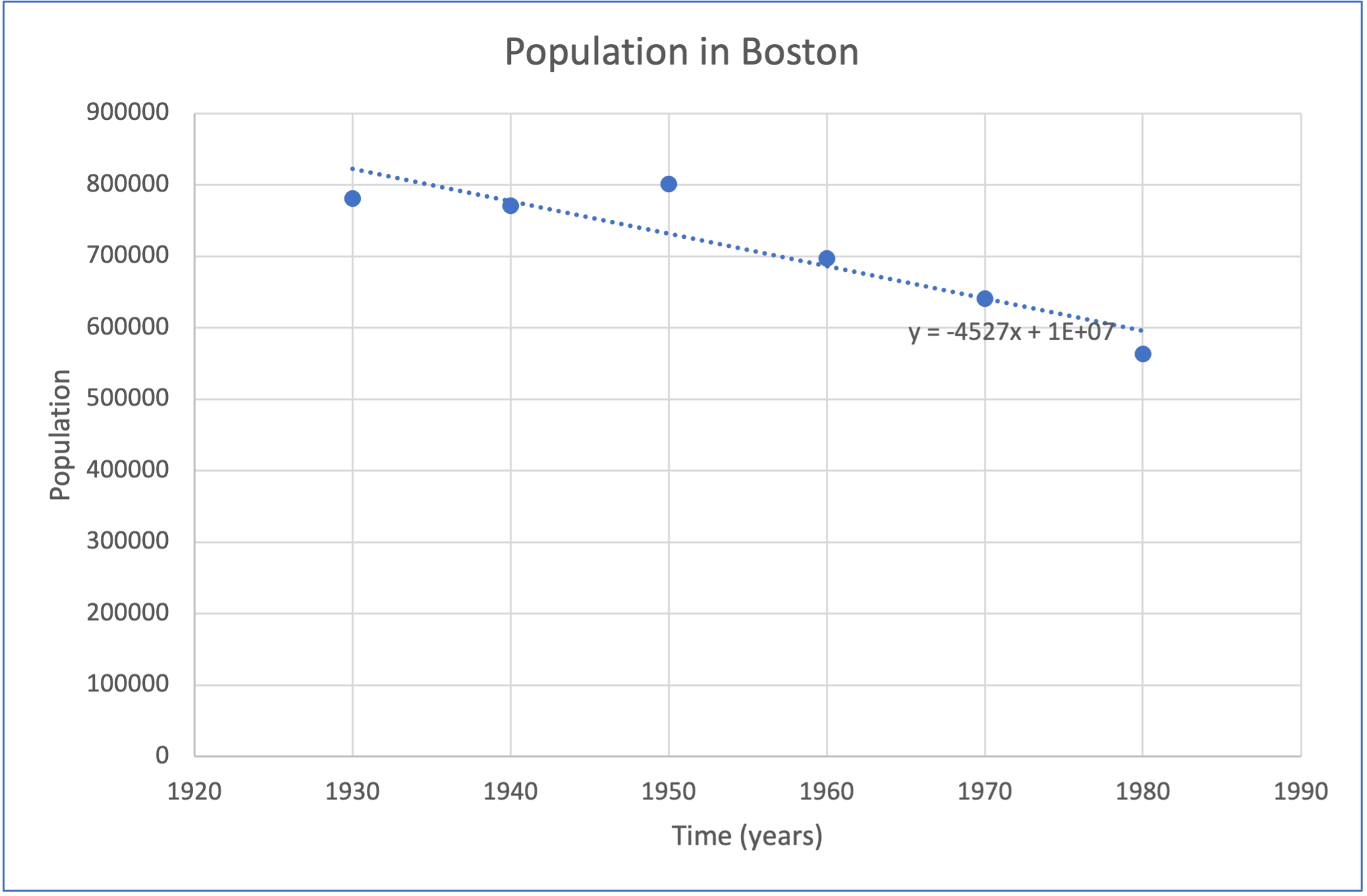}
        \includegraphics[width=0.915\linewidth]{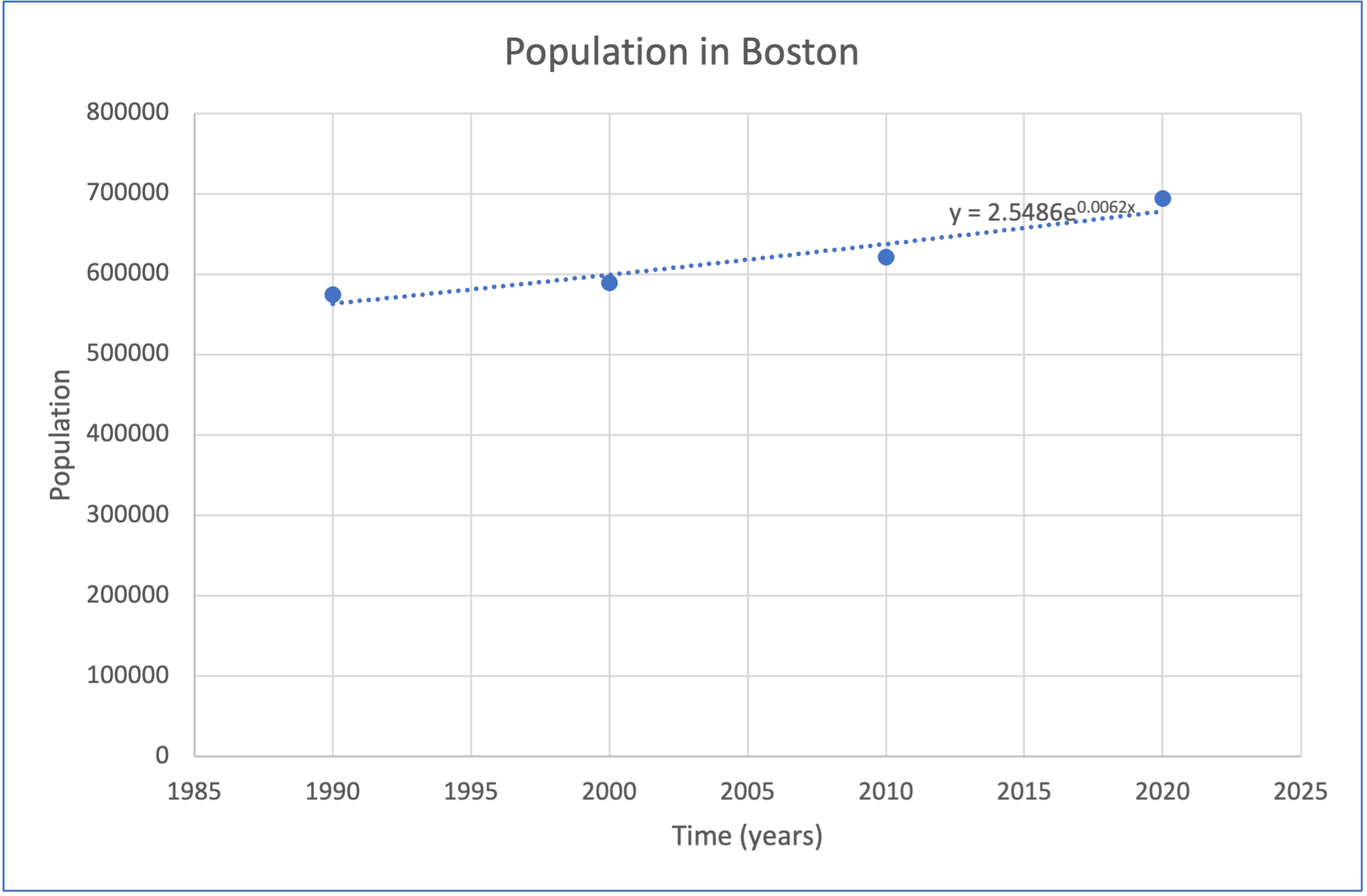}
        \caption{Population of Boston}
     \end{subfigure}
        \caption{Population modeled by exponential functions}
        \label{fig16}
\end{figure}

We use different exponential functions to model them\footnote{See Figure~\ref{fig16} to visualize these functions and see Table~\ref{table3} to know how these functions are defined.}. Since we investigate a wider time range than Meyer does in \cite{12}, we separate the data into five intervals, each of which has a different growth rate. Based on Definition~\ref{def:similar}, we used the exponential model to fit every period, except for period 4 (see Table~\ref{table3}), where the population in Boston and in Manhattan decreased. We use a linear regression line to model this period. 

\begin{table}[!ht]
\begin{center}
\begin{tabular}{ |c|c|c|c|c| }
\hline 
 & \multicolumn{4}{|c|}{Population}\\
\hline 
& \multicolumn{2}{|c|}{Boston}  & \multicolumn{2}{c|}{Manhattan} \\
\hline 
Period & years & equation & equation & years \\
\hline 
1 & 1640-1742 & $P_{B, 1} = 7\cdot10^{-15}\cdot1.0247^t$ & $P_{M, 1} = 7\cdot10^{-12}\cdot1.0203^t$ & 1698-1756 \\
\hline 
2 & 1780-1870 & $P_{B, 2} = 3\cdot10^{-21}\cdot1.0324^t$ & $P_{M, 2} = 5\cdot10^{-27}\cdot1.0405^t$ & 1771-1840 \\
\hline 
3 & 1880-1920 & $P_{B, 3} = 3\cdot10^{-10}\cdot1.0187^t$ & $P_{M, 3} = 7\cdot10^{-14}\cdot1.0239^t$ & 1850-1910 \\
\hline 
4 & 1930-1980 & $P_{B, 4} = -4527\cdot t+10^7$ & $P_{M, 4} = -11889\cdot t+2\cdot10^7$ & 1920-1970 \\
\hline 
5 & 1990-2020 & $P_{B, 5} = 2.5486\cdot1.0062^t$ & $P_{M, 5} = 1941.1\cdot1.0033^t$ & 1980-2019 \\
\hline
\end{tabular}
\caption{The exponential function for modeling the population of the dynamical fractal (Boston $B$) and (Manhattan $M$) of Figure~\ref{fig16}}
\label{table3}
\end{center}
\end{table}

Based on Definition~\ref{def:similar_shape} and Lemma~\ref{lem:avergae_curvature} we calculated all $\alpha$ and compare them (see Table~\ref{tab:alpha}).

\begin{table}[!ht]
\begin{center}
\begin{tabular}{ |c|c|c| }
\hline 
 & $\alpha_{B, i}$ & $\alpha_{M, i}$\\
\hline 
$i = 2$ & 95.9158 & 58.005 \\
\hline 
$i = 3$ & 174.575 & 485.729 \\
\hline 
$i = 5$ & 0.403918 & 0.208647 \\
\hline
\end{tabular}
\caption{The comparison of the $\alpha$ value (defined in Definition~\ref{def:similar_shape}) for Boston and Manhattan (Exceptional case: when $i = 5$, we computed the ratio of the average curvature of period 3 to period 5 because period 4 is modeled by a linear regression line).}
\label{tab:alpha}
\end{center}
\end{table}

When $i = 3$, the curvature of the Manhattan population changed more drastically than that of the Boston population. When $i = 2$ and $i = 5$, they both have similar curvature changes. Also, since in period 4, both the Manhattan population and Boston population started to decrease for 50 years, we believe the long-term trend of the fractal dimension of Manhattan is similar to the long-term trend of Boston's fractal dimension. 

\subsection{Fitting the Data with the Logistic Differential Equation}
\label{sub:atfdb:differential}

In this section, we fit the fractal dimension of Manhattan into the logistic differential equation out of its greater abundance of urban layer images. We use this as an approximation of the long-term trend of the fractal dimension of Boston. 
\begin{figure}[H]
     \centering
     \begin{subfigure}{0.4\textwidth}
         \centering
         \includegraphics[width=0.8\linewidth]{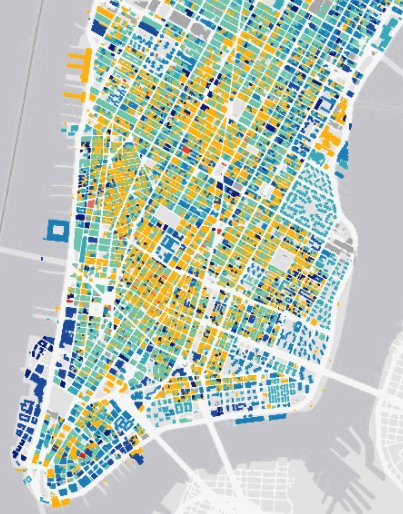}
         \caption{Urban layer of Manhattan in 2000 \cite{26}}
     \end{subfigure}
     \hfill
     \begin{subfigure}{0.49\textwidth}
        \centering
        \includegraphics[width=0.8\linewidth]{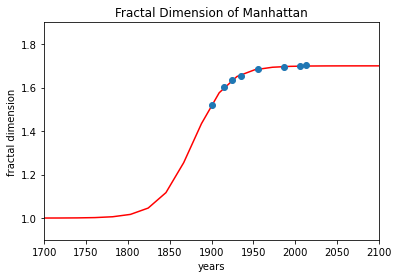}
        \caption{Fractal dimension of Manhattan and the logistic equation}
        \label{17b}
     \end{subfigure}
        \caption{Fitting the fractal dimension by logistic model}
\end{figure}

We tested the urban layers of Manhattan. An urban layer, as an interactive map created by Morphocode \cite{26}, offers a description of Manhattan's urban fabric through the chronological depiction similar to a planimetric map of the spread of buildings while significantly minimizing error. Thus, based on that, we calculated its fractal dimension in 1900, 1915, 1924, 1935, 1956, 1986, 2006, and 2013 and fit it into the logistic differential equation (see Definition~\ref{def:differential}). 
\begin{align}
    \Phi^*_{M, t} = \frac{0.699952}{1+2.07022\cdot10^{40}\cdot e^{-0.049432t}}+1,
    \label{eq:5.1}
\end{align} 
where $M$ represents the dynamical fractal, Manhattan. By plotting $\Phi_M$, we got Figure~\ref{17b}. We set the lower bound of the logistic model as $1$ (as the state before this city was built) and approximate its upper bound as equal to $1.699952$ (as the state where most of the city facilities have been built) using a regression line fit\footnote{The upper bound of the logistic model might vary due to the change of the socioeconomic status, which might be influenced by political, technological, or natural factors \cite{12}.}. Before 1800, the graph shows the Manhattan area was not so urbanized and it was developing slowly. After this period of becoming fully developed and overcrowded, Manhattan reached the carrying capacity of its urbanization described by a gradually stabilizing fractal dimension. 

This model might not always account for social, political, and economic factors, but indicates the trend of and provides an understanding of urbanization. That is because urban development is similar to the spread of the virus, which, logically, should follow the logistic model since the more people get infected, the faster they will infect others, and gradually it will slow down and approach its upper bound. However, in the real world, the spread of the virus hardly follows the logistic model under the influence of the vaccine, governments' policies, variants, and the change of awareness of people, just like COVID-19. Likewise, a positive acceleration of growth might occur because an increase in population and in urbanization might attract more people to the area but might variably be influenced by the government's policies, people's preferences, and resource allocation.

We apply Lemma~\ref{lem:difference_from_differential} to transform Eq.~\eqref{eq:5.1} to difference equation:
\begin{align*}
    x(t+1)\approx 1.044[1-x(t)]x(t).
\end{align*}

Finally, by Lemma~\ref{thm:stable}, since $1.044\in (0,1]$, there is a large neighborhood of stability, and by Theorem~\ref{c:stability}, both the differential and the difference equations are stable, so we expect that in the long run, the fractal dimension of Boston will be convergent to $1.699952$. 

\section{Conclusion}
\label{sec:conclusion}
In conclusion, we have developed the formal definition of dynamical fractal. Since dynamical fractals are changing fractals, we can use the trend of fractal dimension to study these dynamical fractals. The curvature may be used to show that two dynamical fractals are similar. In our case study, it might not be as meaningful to keep analyzing the long-term trend of the fractal dimension of Boston as other developing cities because Boston has already been very urbanized, and developed cities like Boston will soon approach their upper bound (based on the logistic model). 

However, in the future, people can still use the data and methods presented in this paper to analyze the trend and the degree of complexity and urbanization of other cities or some other developing Boston boroughs, both of which are beneficial to Boston's urban planning and public health. Furthermore, future research can keep developing the theoretical concepts of dynamical fractals, apply the concept of dynamical fractals to analyze other real-world fractals, like cells, cancers, and clouds to learn about their long-term trend, and take the height of the building into consideration for analyzing dynamical fractals.

\appendix
\section{Preliminary}
\label{sec:preli}

In Section~\ref{sub:preli:probability}, we introduce the countable set and the probability space and analyze their relationship. In Section~\ref{sub:preli:vector}, we introduce the basic definition related to vectors and curves.

\subsection{Probability}
\label{sub:preli:probability}

In this section, we first introduce the definition of the countable set.

\begin{definition}[Countable set, Definition 5.2.1 in \cite{c16}]\label{def:countable}
    Let $X$ be a set. $X$ is countable if there exists an injective function $f$ such that
    \begin{align*}
        f: X \to \mathbb{N}.
    \end{align*}
    
\end{definition}

\begin{definition}[$\sigma$-algebra, Definition 2.2 in \cite{g05}]\label{def:algebra}
    Let $X$ be a set. The $\sigma$-algebra $\Sigma$ on the set $X$ is a collection of the subsets of $X$, satisfying
    \begin{itemize}
        \item $X \in \Sigma$,
        \item if $A \in \Sigma$, then $A^c \in \Sigma$, and
        \item if for all $i$, $A_i \in \Sigma$, then $\bigcup_i A_i \in \Sigma$.
    \end{itemize}
\end{definition}

\begin{definition}[Probability space, Definition 3.1 in \cite{g05}]\label{def:probability}
    A probability space $(\Omega ,\Sigma,P)$ is a structure, satisfying
    \begin{itemize}
        \item $\Omega$ is the sample space, which is set containing all possible outcomes.
        \item $\Sigma$ is the $\sigma$-algebra (see Definition~\ref{def:algebra}) on $\Omega$ and is the set of events.
        \item $P$ is a function, where $P : \Sigma \to [0, 1]$, satisfying
        \begin{itemize}
            \item For all $A \in \Sigma$, $P(A) \geq 0$,
            \item $P(\Sigma) = 1$, and
            \item For all $A_m, A_n \in \Sigma$, satisfying $A_m \cap A_n = \emptyset$, 
            \begin{align*}
                P(\bigcup_i A_i) = \sum_i P(A_i).
            \end{align*}
        \end{itemize}
    \end{itemize}
\end{definition}

\begin{fact}\label{fac:probability}
    Let $t_0, t_* \in \mathbb{R}$, with $t_0 < t_*$.

    Let $((t_0, t_*), \Sigma_B, P)$ be a probability space (see Definition~\ref{def:probability}), where $\Sigma_B$ is a $\sigma$-algebra satisfying
    \begin{align*}
        \Sigma_B := \{(a, b), [a, b), (a, b], [a, b] \mid a, b \in (t_0, t_*) \},
    \end{align*}
    and $P$ is a probability function defined on $\Sigma_B$, defined as
    \begin{align*}
        P((a, b)) = \frac{b - a}{t_* - t_0}.
    \end{align*}

    Then, for all $t \in (t_0, t_*)$, we have
    \begin{align*}
        P([t, t]) = P(\{t\}) = 0.
    \end{align*}
\end{fact}

\subsection{Vectors and Curves}
\label{sub:preli:vector}

In this section, we present the definitions related to vectors and curves.

\begin{definition}[Unit-speed curve, Definition 1.2.3 in \cite{36}]\label{def:unit_speed}
    If $\gamma : (\alpha, \beta) \to \mathbb{R}^n$ is a parametrized curve, its speed at the point $\gamma(t)$ is $\| \gamma'(t)\|$, and $\gamma$ is said to be a unit-speed curve if $\gamma'(t)$ is a unit vector for all $t \in (\alpha, \beta)$.
\end{definition}

\begin{definition}[Regular point/curve, Definition 1.3.3 in \cite{36}]\label{def:regular}
    A point $\gamma(t)$ of a parametrized curve $\gamma$ is called a regular point if $\dot{\gamma}(t) \neq 0$. A curve is regular if all of its points are regular.
\end{definition}

\begin{definition}[Curvature, Definition 2.1.1 in \cite{36}]\label{def:curvature}
    Let $\gamma$ be a unit-speed curve (see Definition~\ref{def:unit_speed}) with parameter $t$. Its curvature $\kappa(t)$ at the point of $\gamma(t)$ is defined as  $\| \ddot{\gamma}(t) \|$.
\end{definition}

\begin{fact}[Proposition 2.1.2 in \cite{36}]\label{fac:curvature}
    Let $\gamma(t)$ be a regular curve (see Definition~\ref{def:regular}) in $\mathbb{R}^3$. Then, its curvature (see Definition~\ref{def:curvature}) is
    \begin{align*}
        \kappa = \frac{\|\ddot{\gamma} \times \dot{\gamma}\|}{\|\dot{\gamma}\|^3}.
    \end{align*}
\end{fact}

\section*{Acknowledgments}
In memory of Professor Emma Previato, this research is supported by grants from the Undergraduate Research Opportunity Program (UROP) of Boston University. The author sincerely thanks Professor Emma Previato, Professor Dafeng Xu, Alex Novakovic, Min Han, Mu Qiao, Quan Zhou, Xiaohua Zhang, and UROP for their tremendous support, including grants, instructions, suggestions, and coding assistance. Without their significant contributions, this paper would not have been possible.

\bibliographystyle{siamplain}
\bibliography{references}

\end{document}